\RequirePackage{lineno}
\documentclass[10pt]{elsart}
\usepackage{amsmath,amssymb}
\usepackage{enumerate}
\journal{ArXiv}
\newcommand{\p}{^{\prime}}

\newcommand{\smsp}{\phantom{.}}

\newcommand{\lam}{{(\lambda)}}
\newcommand{\mam}{{(\mu)}}

\newcommand{\Z}{\ensuremath{\mathbb{Z}}}
\newcommand{\ma}{\ensuremath{\frak{m}}}

\newcommand{\ess}{\ensuremath{\mathfrak{S}}}
\newcommand{\rmk}{\noindent\textbf{Remark: }}
\newcommand{\at}{\ensuremath{\mathfrak{P}}}
\newcommand{\enn}{{(n)}}
\newcommand{\bk}{\ensuremath\langle k\rangle}
\newcommand{\bv}{\textbf{a}}
\newcommand{\bS}{\mathbb{S}}
\newcommand{\jenn}{{(j,n)}}

\newcommand{\R}{\ensuremath{\mathbb{R}}}
\newcommand{\T}{\ensuremath{\mathbb{T}}}
\newcommand{\N}{\ensuremath{\mathbb{N}}}
\newcommand{\C}{\ensuremath{\mathbb{C}}}
\newcommand{\GG}{\ensuremath{\Gamma(l,n)}}
\newcommand{\DD}{\ensuremath{\Delta(l,n)}}
\newcommand{\pid}{\ensuremath{[-\pi,\pi]^d}}
\newcommand{\BB}{\mathfrak{B}}
\newcommand{\bB}{\mathbb{B}(l,n)}
\newcommand{\bn}{{\mathfrak{B}(n)}}
\newcommand{\bwn}{{\mathfrak{B}(w,n)}}
\newcommand{\ZZ}{\mathfrak{Z}}
\newcommand{\Card}{\text{Card}}
\newtheorem{theorem}{\textbf{Theorem}}[section]
\newtheorem{corollary}[theorem]{\textbf{Corollary}}

\newtheorem{lemma}[theorem]{\textbf{Lemma}}
\numberwithin{equation}{section}
\begin{document}
\begin{frontmatter}
\title{A Central Limit Theorem for Periodograms Under Rho-Prime Mixing}
\author{R. J. Niichel}
\address{Ivy Tech Community College, Bloomington, IN} \ead{rniichel@indiana.edu}
\begin{abstract}
    This paper consists of a proof of a multivariate Central Limit Theorem for ``rectangular'' sums of dependent complex-valued, $\rho\p$-mixing random variables indexed by $\Z^d$.
\end{abstract}
\end{frontmatter}
\section{Introduction}
The fundamental material required for the main results is presented in this section.

\subsection{Historical Introduction}

The Classical Central Limit Theorem states that a normalized sum of independent, identically distributed random variables converges in distribution to a normal random variable.  The result goes back as far as DeMoivre, who did not publish a CLT \emph{per se}, but rather an approximation for the probabilities of normalized sums.  Laplace was the first to publish a result (in 1810, after almost 40 years of work) that could properly be called a CLT.  Over the next 140 years, many new CLTs were published, most with assumptions about the random variables being independent (for an interesting study of the CLT during the 19th Century, see Fischer (2010)).  However, in 1951, Donsker published his weak invariance principle, essentially putting an end to CLT's for independent random variables.

Then, in 1956, Murray Rosenblatt published a paper in which he proved a CLT on a sequence of random variables that were assumed to be dependent.  His method involved developing a way to measure the dependence, and then insisting that if the random variables were far removed from each other (whether in space or time), then their measured dependence should be small.  This, as well as a few additional assumptions, permitted the proof of a CLT.

Since that time, a number of new measures of dependence have been developed, each measuring the dependence in a more or less intuitive way.  This paper deals with what is called the $\rho\p$-mixing condition.

\subsection{The Main concepts}\label{C1S2}
The goal of this section is to introduce the main concepts which will be used throughout the remainder of this paper.

Most of the classical CLTs work on sequences of random variables or triangular arrays.  This notion will be generalized in this paper to include fields of random variables.  A \textbf{{field}} $X$ of random variables is a collection of random variables which are indexed by $\Z^d$ ($d\in\N$); i.e. $X:=\{X_k:k\in\Z^d\}$.  Such a field will be called ``\textbf{complex}'' or ``\textbf{real}'' if the range of the individual random variables is $\C$ or \R, respectively.  If the mean of each random variable is zero (i.e. $EX_k=0$ for all $k$), then the field will be called ``\textbf{centered}''.

Another assumption of the original CLT is that the sequence should be identically distributed.  This is a natural assumption to make in many applications, since one usually assumes that whatever one is measuring behaves the same way on different trials.  In this paper, similar notions will be assumed.  However, since the random variables will be allowed to be dependent, it will be necessary to discuss how they are distributed with respect to each other.  For example, saying that $X$ is ``\textbf{strictly stationary}'' means that for every non-empty subset $S$ of $\Z^d$, and for every $p\in\Z^d$, the subsets $\{X_k:k\in S\}\subset X$ and $\{X_{k+p}:k\in S\}\subset X$ have the same distribution.  In other words, the joint distributions are fixed under translations of the indexing set.  A useful characteristic of strictly stationary fields is that it requires no loss of generality to assume that the field is centered.

A similar criterion for a random field $X$ is weak stationarity.  To be ``\textbf{weakly stationary}'', a field $X$ must satisfy three conditions: 1) for all $k\in\Z^d$, $E|X_k|^2=\sigma^2<\infty$; 2) there exists a $\mu\in\C$ such that $EX_k=\mu$ for all $k\in\Z^d$; 3) $E(X_k-\mu)\overline{(X_j-\mu)}$ depends only on the vector $k-j$.

As it turns out, there is an easy way to generate a centered, complex, weakly stationary (``\textbf{CCWS}'') field from a centered, complex, strictly stationary (``\textbf{CCSS}'') one with finite second moments.  The method is as follows:  First, take a vector $\lambda=(\lambda_1,\lambda_2,\dots,\lambda_d)\in(-\pi,\pi]^d$, and write $e^{i\lambda}=(e^{i\lambda_1},e^{i\lambda_2},\dots,e^{i\lambda_d})$.  (In general, if $k\in\R^d$ is a vector, then $e^{ik}$ will be a shorthand for the vector $(e^{ik_1},e^{ik_2},\dots,e^{ik_d})$).  Then, define the field $X^\lam:=\{X_k^\lam:=e^{-ik\cdot\lambda}X_k:k\in Z^d\}$, where the dot in the exponent is the standard dot product.  It is now easy to see that $X^\lam$ is CCWS; however, $X^\lam$ need not be strictly stationary.

One characteristic of CCWS fields that is of particular interest for this paper is the fact that they often (but not always) have spectral densities.  Let $Y=\{Y_k:k\in\Z^d\}$ be a CCWS field; let $\T$ denote the unit circle in $\C$; and let $\ma_\T:=dz/2\pi i z$ be the normalized Lebesgue measure on $\T$, with $\ma_\T^d=\ma_\T\times\ma_\T\times\dots\times\ma_\T$ denoting the product measure on $\T^d$.  Then, the non-negative Borel function $f:\T^d\to\R$ is the \textbf{spectral density} of $Y$ if for every pair $k\neq j$

\[EY_k\overline{Y}_j=\int_{\T^d}e^{i(k-j)\cdot\theta}f(e^{i\theta})d\ma_\T^d(e^{i\theta}).\]
More will be said about the spectral density in a moment.

The other main concept which needs introduction is the notion of mixing conditions.  Intuitively, many phenomena in the real-world are dependent, and independence is often a lot to ask.  Take, for instance, the disturbance in a radio signal measured every minute.  If there is a large quantity of static in the current signal, our expectation for the next signal is that it will also have a lot of static, perhaps due to a storm or other disturbance.  Thus the time-dependence of this sequence is apparent in small time intervals.  On the other hand, if there is a high level of disruption in the present signal, what can be said about the measurement taken a month or year from now?  In this case, it seems reasonable to assume that measurements separated by large intervals of time should be more or less independent.

Capturing this idea in a rigorous form means first of all that the dependence must be measured.  Consider then the \textbf{(maximum) correlation coefficient} for two sigma-fields $\mathcal{A}$ and $\mathcal{B}$:
\[
\varrho\left(\mathcal{A},\mathcal{B}\right)=\sup_{f,g} Corr(f,g),
\]
where the supremum is taken over all $f\in L^2(\mathcal{A})$ and $g\in L^2(\mathcal{B})$ (real- valued), and $Corr(f,g)=(Efg-EfEg)/\|f\|_2\|g\|_2$.

To apply the correlation coefficient to the context of a random field $X$, take a non-empty subset $V\subset\Z^d$, and let $\sigma(V)$ denote the sigma field generated by the random variables $X_k\in X$ with indices $k\in V$.  Then, define
\[
\varrho\p(X,n)=\sup_{S,T}\varrho\left(\sigma(S),\sigma(T)\right)
\]
where now the supremum is taken over all finite non-empty sets $S$ and $T$ which are separated by $n$ units in (at least) one dimension. That is to say, there is a subscript $u$, $1\leq u\leq d$ so that if $S\ni k=(k_1,\dots,k_d)$ and $T\ni l=(l_1,\dots,l_d)$, then $|k_u-l_u|\geq n$.  It is important to note that the sets $S$ and $T$ can be ``interlaced'', meaning there may be $k,j\in S$ and $l\in T$ such that $k_u\leq l_u\leq j_u$, and/or vice versa.

It is now possible to describe what is meant by ``$\varrho\p$-mixing''.  A field $X$ is said to be \textbf{$\varrho\p$-mixing} if $\varrho\p(X,n)\to 0$ as $n\to\infty$.  Again, what is being said here is that the random variables are ``asymptotically independent'', insofar as $\varrho\p$ measures dependence.

Another, perhaps better-known, measure of dependence is the $\varrho^*$ condition.  It is very similar to $\varrho\p$:
\[
\varrho^*(X,n)=\sup_{S,T}\varrho\left(\sigma(S),\sigma(T)\right).
\]
The only difference is that the elements $k$ and $l$ of the finite non-empty sets $S$ and $T$ (respectively) must satisfy $\|k-l\|\geq n$ (here and below, $\|\cdot\|$ is the standard Euclidean norm):
\[
\min_{k\in S, l\in T}\|k-l\|\geq n.
\]
Note that in one dimension the two mixing conditions are equivalent.

It is easy to see that $\varrho^*$-mixing implies $\varrho\p$-mixing, since $\varrho^*(n)\geq\varrho\p(n)$.  Thus (since only $\varrho\p$-mixing is assumed below) all of the results proved in this paper apply to $\varrho^*$-mixing fields as well.

The final concept which requires introduction is the periodogram.  Let $v=(v_1,v_2,\dots,v_d)\in\N^d$ be an an arbitrary vector.  The $d$-dimensional ``box'' $\BB(v)$ defined by the vector $v$ will be
\[
\frak{B}(v):=\{w\in\N^d:1\leq w_j\leq v_j\smsp, j=1,2,\dots,d\}
\]
Now suppose that $\{v^\enn=(v_1^\enn,v_2^\enn,\dots,v_d^\enn)\}_{n=1}^\infty$ is a sequence of vectors in $\Z^d$. If this sequence is temporarily fixed, then the collection $\{\bn:=\BB(v^\enn)\}_{n=1}^\infty$ is a ``sequence'' of subsets of $\N^d$.  Also, define the sequence of real numbers
\[
V^\enn:=\prod_{j=1}^dv_j^\enn.
\]
Notice that $V^\enn=\Card\left(\bn\right)$.  The \emph{n}th periodogram will sum the random variables whose subscripts lie in the \emph{n}th box.

Next, if $\ess\subset \Z^d$ is nonempty, define
\[
S_\ess^\lam:=\sum_{k\in\ess}X_k^\lam.
\]
Again assuming that $\{v^\enn\}_{n=1}^\infty$ is a fixed sequence, it is possible to define
\[
S_n^\lam:=S_{\bn}^\lam.
\]
Finally (assuming again that $v^\enn$ is a fixed sequence), the \textbf{periodogram} is defined to be
\[
I_n^\lam:=\frac{\left|S_n^\lam\right|^2}{V^\enn}.
\]
As the title suggests, the periodogram is the focus of this paper.  Notice that in this definition, the $n$th periodogram depends \emph{both} on a particular $\lambda$ and on the vector $v^\enn$.  Sometimes the sequence $v^\enn$ is defined in the context of a particular theorem below; other times, it is not given a precise definition.  What should be borne in mind in these latter situations is that the sequence $v^\enn$ is assumed to be arbitrary (up to a stated condition) but fixed.

\subsection{Motivation}

The following is the main result of this paper:

\begin{theorem}\label{N1}
Let $d$ be a positive integer and suppose $X:=\{X_k,k\in\Z^d\}$ is a $\varrho\p$-mixing, CCSS random field such that $E|X_k|^2=\sigma^2<\infty$.  Let $f\lam:=f(e^{i\lambda})$ be the (continuous) spectral density of $X$. Let $\lambda\in\at$, and let $\{\lambda^\jenn\}_{n=1}^\infty$, $j=1,2,\dots,m$ be sequences of elements of \pid which converge to $\lambda$, and which satisfy the conditions of Lemma \ref{Lemma2.5}. Suppose $\{v^\enn\}$ is a sequence of vectors that satisfies \eqref{2.1.min}, i.e.
\[\lim_{n\to\infty}\min\{v_1^\enn,v_2^\enn,\dots,v_d^\enn\}=\infty.
\] Then
\begin{equation}\label{Main3}
    \frac{\bS_n^0}{\sqrt{V^\enn}}\Rightarrow \textbf{Z}
\end{equation}
where $\textbf{Z}:\Omega\rightarrow \R^{2m}$ has the normal distribution with the $2m\times2m$ covariance matrix
\[\Upsilon_m^\lam:=
\begin{bmatrix}
f\lam&0&\dots&0\\
0&f\lam&\dots&0\\
\vdots&\vdots&\ddots&\vdots\\
0&0&\dots&f\lam
\end{bmatrix}.
\]
\end{theorem}

The argument is fairly standard, and in fact very similar to a result which the author has already published (see \cite{13}).  Theorem \ref{N1} makes it possible to prove (in a forthcoming paper) a weak-type law for periodograms, which does not require a summable cumulant assumption.  The purpose of presenting it on arXiv is to avoid having to give the long argument for Theorem \ref{N1} in a formal paper.

The main motivation behind proving Theorem \ref{N1} was to remove some of the restrictions of a Central Limit Theorem proved by Murray Rosenblatt (1985).  There, Rosenblatt proved a CLT for $\rho^*$-mixing fields ($\rho\p$-mixing was unknown at the time), relying on finite summable second and fourth order cumulants to prove his result. (For a definition of cumulants, please see Rosenblatt (1985), page 33ff.)  The reader will notice that Theorem \ref{N1} does not rely on cumulants, only the distance between the sampling frequencies (the $\lambda^{(j,n)}$'s) to achieve the CLT.

\subsection{Important Results}\label{C1S4}
There are a number of important results which pertain to the subject matter of this paper.  Some of them are well-known, others perhaps not so much.

Start with the familiar:

\begin{lemma}[Slutsky]\label{slutsky}
Let, $Y$, $Y_n$, $Z_n$, $n=1,2,\dots$, be random variables, and let $\{a_n\}_{n=1}^\infty$ be a sequence of real numbers.  Suppose that $Y_n\Rightarrow Y$, $a_n\to a$, and $Z_n\Rightarrow0$ as $n\to\infty$.  Then $a_nY_n+Z_n\Rightarrow aY$ as $n\to\infty$.
\end{lemma}

The next theorem is crucial to the workings of the Bernstein blocking argument (see \S\ref{C3S4}); it will permit a reduction of the main CLT which deals with dependent random variables to a CLT that deals with independent random variables.

\begin{theorem}[Billingsley, Theorem 26.3]\label{Bill26.3}  Let $\mu_n,\mu$ be probability
measures on $(\mathbb{R},\mathcal{R})$ ($\mathcal{R}$ is the Borel
sigma-field on $\mathbb{R}$) with characteristic functions (Fourier
transforms) $\phi_n$ and $\phi$. A necessary and sufficient
condition for $\mu_n\Rightarrow\mu$ is that $\phi_n(t)\to\phi(t)$
for every $t$.
\end{theorem}

Also of fundamental importance to the blocking argument is the Cramer-Wold device, which basically turns a multivariate CLT into a univariate CLT:

\begin{theorem}[The Cramer-Wold Device]\label{Cramer} For
$\mathbb{R}^k$-valued random vectors $X^{(n)}=(X_1^{(n)},\hdots
X_k^{(n)})$ and $Y=(Y_1,\hdots Y_k)$, a necessary and sufficient
condition for $X^{(n)}\Rightarrow Y$ is that
$\sum_{u=1}^kt_uX_u^{(n)}\Rightarrow\sum_{u=1}^kt_uY_u$ for each
$(t_1,\hdots t_k)$ in $\mathbb{R}^k$. \end{theorem}

The Mapping Theorem is handy once a CLT has been proved:

\begin{theorem}[The Mapping Theorem, \cite{1}, Thm. 29.2]\label{mapp} Suppose that $h:\R^k\to\R^j$ is a measurable function, and suppose that the convergence of measures $P_n\Rightarrow P$ holds in $\R^k$.  Let $D_h$ be the set of discontinuities of $h$.  If $P(D_h)=0$, then $P_n h^{-1}\Rightarrow P h^{-1}$.\end{theorem}

Another way to state Theorem \ref{mapp} is to say that if $X_n\Rightarrow X$, and $P_n$ and $P$ are their respective distributions, and if $P(D_h)=0$, then $h(X_n)\Rightarrow h(X)$ (see \cite{1.75}).

To prove Theorem \ref{N1}, the following result will be essential:

\begin{theorem}[Billingsley, Theorem 25.12]\label{Bill25.12} Let $r$ be a positive
integer.  If $Y_n\Rightarrow Y$ and
$\sup_n\{E|Y_n|^{r+\epsilon}\}<\infty$ for some $\epsilon>0$, then
$E|Y|^r<\infty$ and $E|Y_n|^r\to E|Y|^r$.
\end{theorem}

The remaining results are not as well-known but nevertheless are quite useful in the study of mixing conditions.  The first result provides a useful means of controlling the second moment of a sum of random variables when the field is $\varrho\p$-mixing.  It is of a type of inequalities known by the moniker ``Rosenthal''.  It is very useful in connection with Lyapounov's criterion.

\begin{theorem}[Bradley, Theorem 29.30]
\label{Brad2}
Suppose $\beta$ is a number in the interval $[2,\infty)$.
Suppose further that $X$ is a (not necessarily stationary) random field of complex-valued random variables
such that for each $k\in\Z^d$, $EX_k=0$, and $E|X_k|^\beta<\infty$.  Suppose that $\varrho\p\enn<1$ for some $n\in\Z$.
Then, for any finite set $\mathfrak{S}\subset\Z^d$, it holds that
\[E\left|\sum_\mathfrak{S} X_k\right|^q\leq C\cdot \left[\sum_\mathfrak{S}
E|X_k|^\beta+\left(\sum_\mathfrak{S} E|X_k|^2\right)^{\beta/2}\right],\]
where C is a constant that depends on d, n, $\varrho\p(n)$, and $\beta$.
\end{theorem}

The final major result may provide some intuition about why the periodogram behaves the way it does asymptotically.  Its proof resembles the standard Fejer Theorem (see \cite{9}).

\begin{theorem}[Bradley, \cite{2}, Theorem 28.21]
\label{Brad1992} Let $v^\enn\equiv(n,n,\dots,n)$. If $X$ is a CCWS random field, such that $\varrho\p(n)\to0$
as $n\to\infty$, then $X$ has a \emph{continuous} spectral density
$f\lam$ on $(-\pi,\pi]^d$ and
\[\lim_{n\to\infty}
EI_n^\lam=f\lam\] and the
convergence is uniform over all $\lambda\in(-\pi,\pi]^d$
\end{theorem}\vspace{.2in}

The following is a minor result needed only for reference:

\begin{lemma}
\label{LemmaZ} Suppose that $a\in[0,\infty)$ and $\{a_k:k\in\Z^d\}$
is a field of non-negative numbers such that for every $\epsilon>0$
there exists an $M=M(\epsilon)>0$ so that whenever the Euclidean norm of
any vector $k$ is greater than $M$,  $|a_k-a|<\epsilon$. Then,
as $n\to\infty$,
\[\frac{\sum_{k\in\bn}a_k}{V^\enn}\to a.\]
\end{lemma}

\section{Preliminary Groundwork}\label{C2}
The first order of business is to prove a few results about the moments of the normalized sums.  For the most part, this involves reproving theorems proved by others which deal only with the special sequence of vectors $v^\enn=(n,n,\dots,n)$ to a more general sequence which only satisfies
\begin{equation}\label{2.1.min}
\lim_{n\to\infty}\min\{v_1^\enn,v_2^\enn,\dots,v_d^\enn\}=\infty.
\end{equation}
This won't be very difficult, but it does need to be done.
\subsection{Reworking Theorem \ref{Brad1992}}\label{C2S3}
The first theorem which needs adjusting is the second half Theorem \ref{Brad1992}:

\begin{theorem}
\label{S1T1}
Let $X:=\{X_k:k\in\Z^d\}$ be a CCWS field of random variables.  Suppose that $\varrho\p(X,n)\to0$ as $n\to\infty$ and let $f\lam:=f(e^{i\lambda})$ denote the (continuous) spectral density.  Suppose $\{v^\enn\}_{n=1}^\infty$ is a sequence of vectors from $\N^d$ which satisfy \eqref{2.1.min}.
\[
\lim_{n\to\infty}I_n^\lam=f\lam
\]
for all $\lambda\in(-\pi,\pi]^d$.  Moreover, this convergence is uniform over all $\lambda$.
\end{theorem}

\pf The proof consists of a simple adaptation of the proof of the standard Fejer Theorem (see \cite{9}, p.176).

Consider $E|S_n^\lam|^2$ and use the spectral density:
\begin{equation}
\label{eq2.1.2}
\begin{aligned}
E|S_n^\lam|^2
&=E\left[\left(\sum_{j\in\bn}e^{-i(j\cdot \lambda)}X_j\right)\left(\sum_{k\in\bn}e^{i(k\cdot \lambda)}\overline{X}_k\right)\right]\\
&=\sum_{j\in\bn}\sum_{k\in\bn}\left[\int_{\T^d}e^{i(j-k)\cdot(\theta-\lambda)}f(\theta)d\ma_\T^d(e^{i\theta})\right]\\
&=\int_{\T^d}\left[\sum_{j\in\bn}\sum_{k\in\bn}e^{i(j-k)\cdot(\theta)}\right]f(\theta+\lambda)d\ma_\T^d(e^{i\theta}).
\end{aligned}
\end{equation}
The sum in the integrand can be transformed using standard algebra:
\begin{equation}
\label{eq2.1.3}
\begin{aligned}
&\sum_{j\in\bn}\sum_{k\in\bn}e^{-i(j-k)\cdot(\theta)}\\
=&\prod_{s=1}^d\left(\sum_{j(s)=1}^{v(n,s)}\sum_{k(s)=1}^{v(n,s)}\exp\{-ij_s\theta_s\}\exp\{ik_s\theta_s\}\right)\\
=&\prod_{s=1}^d\left(\left|\frac{1-e^{iv(n,s)\theta}}{1-e^{i\theta}}\right|^2\right).
\end{aligned}
\end{equation}
(The notations $j(s)=j_s$, $k(s)=k_s$, $v(n,s)=v_s^\enn$ have been substituted here (and elsewhere) to avoid embedded subscripts and superscripts, which tend to be too small to read).

Now is an opportune moment to introduce the Fejer Kernels.  A Fejer Kernel $K(\alpha,n):\R\to\R$ is defined to be
\begin{equation}
\label{F1}
K(\alpha,n):=\frac{\sin^2(n\alpha/2)}{n\sin^2(\alpha/2)}.
\end{equation}
The Fejer Kernels have a number of convenient properties.  Among them are the following (see \cite{11}):

\begin{equation}
\label{F3}
K(\alpha,n)\leq\frac{\pi}{n\alpha^2},\quad 0<|\alpha|<\pi
\end{equation}

\begin{equation}
\label{F4}
\int_{-\pi}^\pi K(\alpha,n)d\alpha=2\pi\quad\text{for all }n
\end{equation}
And, the important thing for this proof is that
\begin{equation}
\label{F5}
\left|\frac{1-e^{in\alpha}}{1-e^{i\alpha}}\right|^2=K(\alpha,n)
\end{equation}
(see \cite{11}, p.30 or \cite{1.75}, Lemma 8.18).

If equations \eqref{eq2.1.2} and \eqref{eq2.1.3} are combined, the following equation holds:
\begin{equation}
\label{eq2.1.8}
EI_n^\lam=\int_{\T^d}\left[\prod_{s=1}^dK\left(\theta_s,v_s^\enn\right)\right]f(\theta+\lambda)d\ma_\T^d(e^{i\theta})
\end{equation}
It must be proved that the integral in \eqref{eq2.1.8} converges to $f\lam$.  So, choose $1>\epsilon>0$.  The function $f$ is bounded on $\T^d$, so $f\lam\leq M$ for all $\lambda$.  Let $|\cdot|_\T$ denote the usual norm on the $d$-torus.  The spectral density is also uniformly continuous on $\T^d$, so let $N_0\in\N$ be such that if $|x-y|_\T<1/\sqrt[3]{N_0}$, then $|f(x)-f(y)|<\epsilon$.  Finally, since \eqref{2.1.min} holds for the sequence $\{v^\enn,\}$, let $N_1\in\N$ be such that

(i) for all $n\geq N_1$, $\sqrt{(v_1^\enn)^{-2/3}+\dots+(v_d^\enn)^{-2/3}}\leq\frac1{N_0}$

(ii) $2M/\sqrt[3]{V^\enn}\leq \epsilon$.

Let the domain of integration in \eqref{eq2.1.8} now be broken into the following two sets:
\[
\begin{aligned}
B(n):=&\left(\frac{-1}{\sqrt[3]{v_1^\enn}}\smsp,\frac{1}{\sqrt[3]{v_1^\enn}}\right)\times\dots\times\left(\frac{-1}{\sqrt[3]{v_d^\enn}}\smsp,\frac{1}{\sqrt[3]{v_d^\enn}}\right)\\
G(n):=&\T^d\setminus B\enn.
\end{aligned}
\]
(Notice that the condition in (i) above ensures that the corners of $B\enn$ are not too far away from the origin.)  Then, using \eqref{F3} and \eqref{F4} with \eqref{eq2.1.8}, together with the definitions of each set,
\[\begin{aligned}
&\int_{B\enn}\left[\prod_{s=1}^dK\left(\theta_s,v_s^\enn\right)\right]|f(\theta+\lambda)-f\lam|d\ma_\T^d(e^{i\theta})\\\leq &\epsilon\int_{B\enn}\prod_{s=1}^dK\left(\theta_s,v_s^\enn\right)d\ma_\T^d(e^{i\theta})
\leq\epsilon
\end{aligned}\]
and
\[
\int_{G\enn}\left[\prod_{s=1}^dK\left(\theta_s,v_s^\enn\right)\right]|f(\theta+\lambda)-f\lam|d\ma_\T^d(e^{i\theta})
<\frac{2M}{\sqrt[3]{V^\enn}}\leq \epsilon.
\]
Hence,
\[
\begin{aligned}
|EI_n^\lam-f\lam|&\leq 2\epsilon.
\end{aligned}
\]
which proves that \eqref{eq2.1.8} converges to $f\lam$.  Notice, however, that this bound is independent of the choice of $\lambda$ since $f$ is uniformly continuous on the domain $\T^d$.  Thus, the convergence is uniform for all $\lambda\in\pid$
\qed

\subsection{The Limit of the Product of Sums}\label{C2S1}

This and the next three sections deal with various lemmas which will be needed to prove the Lemma of \S\ref{C2S6}, which itself will be needed for the proof of the main CLT of Section \ref{C3}.  Throughout the rest of the paper, define the set
\[\at:=\pid\setminus\{-\pi,0,\pi\}^d.\]

\begin{lemma}
\label{lemma2.2}
Let $X$ be a $\varrho\p$-mixing random field such that $EX_kX_j=EX_{k-j}X_0$.  Suppose that $\lambda\in\at$, and let $\{\lambda^{(1,n)}\}_{n=1}^\infty$ and $\{\lambda^{(2,n)}\}_{n=1}^\infty$ be two sequences which converge to $\lambda$.  If the sequence $\{v^\enn\}$ satisfies \eqref{2.1.min}, then
\begin{equation}
\label{lemma2.2.1}
\lim_{n\to\infty}E\left[\frac{S_n^{\lambda(1,n)}S_n^{\lambda(2,n)}}{V^\enn}\right]=0
\end{equation}
\end{lemma}

\pf Without loss of generality, suppose that $e^{i\lambda_1}$ is the coordinate of $e^{i\lambda}$ which is neither 1 nor -1.  Let $e_1$ denote the unit vector (1,0,0,\dots,0). Define the following subsets of $\Z^d$:
\[
\begin{aligned}
A(n)&:=\{j-e_1:j\in\bn\}\\
B(n)&:=\bn\cap A\enn\\
C(n)&:=\bn\setminus B\enn\\
D(n)&:=A\enn\setminus B\enn
\end{aligned}
\]
(Note that $B\enn$ above is unrelated to $B\enn$ from the proof of Theorem \ref{S1T1}.)  Notice that
$\Card(C\enn)=\Card(D\enn)=\prod_{j=2}^d v_j^\enn$; this number will be denoted by $\Lambda^\enn$.  It is also easy to see that  $\Card(B\enn)=(v_1^\enn-1)\Lambda^\enn$.

With the above notations, it is trivial to show that
\begin{equation}
\label{eq2.2.2}
\begin{aligned}
&S_{\bn}^{\lambda(1,n)}S_{\bn}^{\lambda(2,n)}-S_{A\enn}^{\lambda(1,n)}S_{A\enn}^{\lambda(2,n)}\\&\quad=
S_{C\enn}^{\lambda(1,n)}S_{C\enn}^{\lambda(2,n)}+S_{C\enn}^{\lambda(1,n)}S_{B\enn}^{\lambda(2,n)}+
S_{B\enn}^{\lambda(1,n)}S_{C\enn}^{\lambda(2,n)}\\&\quad\quad-S_{D\enn}^{\lambda(1,n)}S_{D\enn}^{\lambda(2,n)}-
S_{D\enn}^{\lambda(1,n)}S_{B\enn}^{\lambda(2,n)}-S_{B\enn}^{\lambda(1,n)}S_{D\enn}^{\lambda(2,n)}
\end{aligned}
\end{equation}
Equation \eqref{eq2.2.2}, together with the Cauchy-Schwarz inequality and Theorem \ref{Brad2} (take $\beta=2$), implies that
\begin{equation}
\label{eq2.2.3}
\begin{aligned}
\left|ES_{\bn}^{\lambda(1,n)}S_{\bn}^{\lambda(2,n)}-ES_{A\enn}^{\lambda(1,n)}S_{A\enn}^{\lambda(2,n)}\right|
&\leq K\|X_0\|_2^2\left[2\Lambda^\enn+4\Lambda^\enn\sqrt{v_1^\enn-1}\right]\\
&<6K\Lambda^\enn\sqrt{v_1^\enn}\|X_0\|_2^2.
\end{aligned}
\end{equation}
Now notice that because $EX_kX_j=EX_{k-j}X_0$, \[ES_{A\enn}^{\lambda(1,n)}S_{A\enn}^{\lambda(2,n)}=\exp\{i\lambda_1^{(1,n)}\}\exp\{i\lambda_1^{(2,n)}\}ES_{\bn}^{\lambda(1,n)}S_{\bn}^{\lambda(2,n)}.\]
This implies that \eqref{eq2.2.3} can be rewritten
\[
\begin{aligned}
&\left|ES_{\bn}^{\lambda(1,n)}S_{\bn}^{\lambda(2,n)}-ES_{A\enn}^{\lambda(1,n)}S_{A\enn}^{\lambda(2,n)}\right|\\
=&\left|1-\exp\{i\lambda_1^{(1,n)}\}\exp\{i\lambda_1^{(2,n)}\}\right|\left|ES_{\bn}^{\lambda(1,n)}S_{\bn}^{\lambda(2,n)}\right|\\
<&6K\Lambda^\enn\sqrt{v_1^\enn}\|X_0\|_2^2,
\end{aligned}
\]
and since $\exp\{i\lambda_1^{(1,n)}\}\exp\{i\lambda_1^{(2,n)}\}\to e^{2i\lambda_1}\neq 1$, it is possible to solve for the expectation.  This proves the theorem, since
\begin{equation}
\label{eq2.2.4}
\left|ES_{\bn}^{\lambda(1,n)}S_{\bn}^{\lambda(2,n)}\right|<
\frac{6K\Lambda^\enn\sqrt{v_1^\enn}\|X_0\|_2^2}{\left|1-\exp\{i\lambda_1^{(1,n)}\}\exp\{i\lambda_1^{(2,n)}\}\right|}
=o\left(V^\enn\right).
\end{equation}
 \qed
\subsection{An Extension}\label{C2S2}
It is not difficult to see that there is a bit more that could be done to improve Lemma \ref{lemma2.2.1} without a significant amount of effort.  In particular, notice that the proof above required that $1-\exp\{i\lambda_s^{(1,n)}\}\exp\{i\lambda_s^{(2,n)}\}\neq0$ for some $s$, $1\leq s\leq d$.  Therefore, the lemma should extend to sets (as opposed to sequences) over which any one of the functions $|1-x_sy_s|$ ($x,y\in\T^d$) has a lower bound.

Recall therefore that
\[\at:=\pid\setminus\{-\pi,0,\pi\}^d.\]
If $K_1$ and $K_2$ are compact subsets of \at\smsp (they must be compact when considered as subsets of $\R^d$), then at least one of the functions
\[g_s:\at\times\at\to[0,2]\quad 1\leq s\leq d\]
\[g_s(\lambda,\mu):=\left|1-\exp\{-i\lambda_s\}\exp\{-i\mu_s\}\right|\]
is bounded below by some $\delta>0$ on $K_1\times K_2$.  Therefore, it is possible to make the bound in \eqref{eq2.2.4} hold for arbitrary $\lambda$ and $\mu$ in $K_1$ and $K_2$.
\[
|ES_n^\lam S_n^{(\mu)}|<
\frac{6K\Lambda^\enn\sqrt{v_1^\enn}\|X_0\|_2^2}{\left|1-\exp\{-i\lambda_s\}\exp\{-i\mu_s\}\right|}
<\frac{6K\Lambda^\enn\sqrt{v_1^\enn}\|X_0\|_2^2}{\delta}.
\]
Moreover, this bound is \emph{uniform} over all $(\lambda,\mu)\in K_1\times K_2$. This is most of the proof of the following corollary:
\begin{lemma}
\label{lemma2.3}
Let $X$ be a $\varrho\p$-mixing random field which satisfies $EX_kX_j=EX_{k-j}X_0$. Let $K_1, K_2$ be subsets of $\at$ which are compact subsets of $\R^d$.  Also let $v^\enn$ be a sequence that satisfies \eqref{2.1.min}.  Then,
\begin{enumerate}[\upshape (i)]
\item The function
\[
f_n\lam:=E\left[\frac{\left(S_n^\lam\right)^2}{V^\enn}\right]
\]
converges uniformly to zero over $K_1$ as $n\to\infty$.\\\\

\item The function
\[
F_n(\lambda,\mu):=E\left[\frac{S_n^\lam S_n^{(\mu)}}{V^\enn}\right]
\]
converges uniformly to zero over $K_1\times K_2$.\\\\

\item Moreover, if $w\in\Z^d$ is a fixed vector, let $\bwn:=\{k+w:k\in\bn\}$. With this notation, the functions
\[
F_{w,n}(\lambda,\mu):=E\left[\frac{S_\bwn^\lam S_\bwn^{(\mu)}}{V^\enn}\right]
\]
converge uniformly to zero at the same rate as $F_n(\lambda,\mu)$; i.e. $|F_n(\lambda,\mu)|=|F_{w,n}(\lambda,\mu)|$.
\end{enumerate}
\end{lemma}

\rmk The statement in (iii) seems to come out of nowhere.  However, it will be very useful to have in the proof of the CLT of Chapter \ref{C3}.

\pf It is easy to see that (i) is a special case of (ii).  Also, (ii) was justified in the work that preceded the statement of the lemma.  The only item which must be proved is (iii).

To do so, notice that $X_{k+w}^\lam=e^{iw\cdot\lambda}\left(e^{ik\cdot\lambda}X_{k+w}\right)$.  And since by assumption $EX_{k+w}X_{j+w}=EX_kX_j$, it follows that
\begin{equation}\label{eq2.3.b} E\left[e^{iw\cdot\lambda}e^{iw\cdot\mu}S_\bwn^\lam S_\bwn^{(\mu)}\right]=e^{iw(\lambda+\mu)}E\left[S_n^\lam S_n^{(\mu)}\right].\end{equation}
The result (iii) follows.\qed

Let $I_\bwn^\lam=\left|S_\bwn^\lam\right|^2/V^\enn$.  The following claim follows for (basically) the same reason that \eqref{eq2.3.b} holds:
\begin{claim}\label{Claim0}If $X$ is CCWS, then for all $\lambda\in\pid$, $E I_\bwn^\lam=EI_n^\lam.$
\end{claim}

\subsection{The Limit of the Covariance of Sums}\label{C2S4}
The previous two sections dealt with the product of $S_n^\lam$ and $S_n^{(\mu)}$.  In the next two sections, the more difficult problem of the covariances is taken up.  I.e., what can be said about the end behavior of \[E\left[\frac{S_n^\lam \overline{S_n^{(\mu)}}}{V^\enn}\right],\] especially if $\lambda$ and $\mu$ are again allowed to vary?  The initial answer to this question is the following lemma:

\begin{lemma}
\label{lemma2.3}
Suppose that $X$ is a CCWS random field which is $\varrho\p$-mixing and that $\{v^\enn\}$ is a sequence which satisfies \eqref{2.1.min}.  Let $\{\lambda^\enn\}_{n=1}^\infty$ and $\{\mu^\enn\}_{n=1}^\infty$ be sequences of elements of $[-\pi,\pi]^d$ such that there are a $\delta$, $0<\delta<.5$, and an $N\in\N$ so that whenever $n\geq N$, there is at least one subscript $s$ so that
\begin{equation}
\label{eq2.4.1}
\left|\lambda_s^\enn-\mu_s^\enn\right|>\frac1{\left(v_s^\enn\right)^{1/2-\delta}}.
\end{equation}
Then,
\[\lim_{n\to\infty}E\left[\frac{S_n^{(\lambda\enn)} \overline{S_n^{(\mu\enn)}}}{V^\enn}\right]=0.\]
\end{lemma}
\pf Choose $\epsilon>0$.  Let $f(e^{i\lambda})=f\lam$ denote the continuous spectral density of the field $X$.  Since $f$ is continuous on \pid, $f(\cdot)\leq M$ for some $M\in \R^+$.  Therefore, choose $n$ so that
\begin{equation}
\label{eq2.4.2}
\frac{3M\sqrt\pi}{(v_j^\enn)^\delta}\leq\epsilon\quad\text{for all }j=1,2,\dots,d.
\end{equation}

Suppose momentarily that $\lambda$ and $\mu$ are any two elements of \pid.  Like in equation \eqref{eq2.1.2},
\begin{equation}
\label{eq2.4.3}
\begin{aligned}
&E\left[S_n^\lam \overline{S_n^\mam}\right]\\
=&\int_{\T^d}\left[\sum_{j\in\bn}\sum_{k\in\bn}\exp\{-i(j\cdot(\lambda-\theta)-k\cdot(\mu-\theta))\}\right]
f(\theta)d\ma_\T^d(e^{i\theta}).
\end{aligned}
\end{equation}
The summation in the integrand can be simplified again:
\begin{equation}
\label{eq2.4.4}
\begin{aligned}
&\sum_{j\in\bn}\sum_{k\in\bn}\exp\{-i(j\cdot(\lambda-\theta)-k\cdot(\mu-\theta))\}\\
=&\prod_{s=1}^d\left[\exp\{-i(\lambda_s-\theta_s)\}\frac{1-\exp\{-iv_s^\enn(\lambda_s-\theta_s)\}}{1-\exp\{-i(\lambda_s-\theta_s)\}}\right.\\
&\quad\quad\quad\left.\times\quad\exp\{i(\mu_s-\theta_s)\}\frac{1-\exp\{iv_s^\enn(\mu_s-\theta_s)\}}{1-\exp\{i(\mu_s-\theta_s)\}}\right]
\end{aligned}
\end{equation}
Letting
\begin{equation}
\label{DKern}
D(\alpha,n):=\frac1{\sqrt n}\cdot\frac{1-e^{-in\alpha}}{1-e^{-i\alpha}},
\end{equation}
it is now clear that
\begin{equation}
\label{eq2.4.6}
\begin{aligned}
&E\left[\frac{S_n^{\lambda\enn} \overline{S_n^{\mu\enn}}}{V^\enn}\right]\\
=&C\int_{\T^d}\left[\prod_{s=1}^dD\left(\lambda_s^\enn-\theta_s,v_s^\enn\right)D\left(\theta_s-\mu_s^\enn,v_s^\enn\right)\right]f(\theta)d\ma_\T^d(e^{i\theta})
\end{aligned}
\end{equation}
where $|C|=1$.

To deal with the integral in \eqref{eq2.4.6}, first assume without loss of generality that the first coordinates of $\lambda$ and $\mu$ satisfy \eqref{eq2.4.1}.  Since $\T^d$ is compact, Fubini's theorem applies, and so the first task will be to control
\begin{equation}
\label{eq2.4.7}
\begin{aligned}
&\left|\int_\T D(\lambda_1^\enn-\theta_1,v_1^\enn)D(\theta_1-\mu_1^\enn,v_1^\enn)f(\theta)d\ma_\T(e^{i\theta})\right|\\
\leq&\int_\T \left|D(\lambda_1^\enn-\theta_1,v_1^\enn)\right|\cdot\left|D(\theta_1-\mu_1^\enn,v_1^\enn)\right|f(\theta)d\ma_\T(e^{i\theta})
\end{aligned}
\end{equation}
(recall that the spectral density is non-negative).  Let $B(x,r)$ denote the ball (in \R) with center $x$ and radius $r$.  Then define the following sets:
\begin{equation}
\label{eq248}
\begin{aligned}
B_1^\enn:&=B\left(\lambda_1^\enn,\frac1{3(v_1^\enn)^{1/2-\delta}}\right)\\
B_2^\enn:&=B\left(\mu_1^\enn,\frac1{3(v_1^\enn)^{1/2-\delta}}\right)\\
B_0^\enn:&=[-\pi,\pi]\setminus(B_1^\enn\cup B_2^\enn).
\end{aligned}
\end{equation}
If you now compare \eqref{F5} and \eqref{DKern}, it is not hard to see that
\[
\left|D(\alpha,n)\right|^2=K(\alpha,n)
\]
(again, $K$ is the Fejer Kernel).  Therefore, the following bounds hold (\emph{cf.} \eqref{F3} and \eqref{eq2.4.1}):
\begin{equation}
\label{eq2.4.8}
\begin{aligned}
\left|D(\lambda_1^\enn-\theta_1,v_1^\enn)\right|&\leq\frac{3\sqrt\pi(v_1^\enn)^{1/2-\delta}}{\sqrt{v_1^\enn}}=\frac{3\sqrt\pi}{(v_1^\enn)^\delta}\quad \text{for } \theta_1\in[-\pi,\pi]\setminus B_1^\enn\\
\left|D(\mu_1^\enn-\theta_1,v_1^\enn)\right|&\leq\frac{3\sqrt\pi}{(v_1^\enn)^\delta}\quad \text{for } \theta_1\in[-\pi,\pi]\setminus B_2^\enn
\end{aligned}
\end{equation}
These bounds suggest that the domain of integration in \eqref{eq2.4.7} should be broken into three pieces, $B_1,^\enn$, $B_2^\enn$ and $B_0^\enn$.  Doing so yields the following estimate
\begin{equation}
\label{eq2.4.9}
\begin{aligned}
&\int_\T \left|D(\lambda_1^\enn-\theta_1,v_1^\enn)\right|\cdot\left|D(\mu_1^\enn-\theta_1,v_1^\enn)\right|f(\theta)d|\ma_\T|(e^{i\theta_1})\\
\leq&
\left[\frac{3M\sqrt\pi}{(v_1^\enn)^\delta}\int_{B_1^\enn}\left|D(\lambda_1^\enn-\theta_1,v_1^\enn)\right|d|\ma_\T|(e^{i\theta_1})\right.\\
&\quad+\frac{3M\sqrt\pi}{(v_1^\enn)^\delta}\int_{B_2^\enn}\left|D(\theta_1-\mu_1^\enn,v_1^\enn)\right|d|\ma_\T|(e^{i\theta_1})\\
&\quad\quad+\left.\frac{9M^2\pi}{(v_1^\enn)^{2\delta}}\int_{B_0^\enn}d|\ma_\T|(e^{i\theta_1})\right]\\
\leq&
2\epsilon\sqrt{\int_\T K(\theta_1,v_1^\enn)d\ma_\T(e^{i\theta_1})}+\epsilon^2\\
\leq&3\epsilon
\end{aligned}
\end{equation}
(the penultimate inequality follows from the Cauchy-Schwarz inequality).

With the help of Cauchy-Schwarz (again), the remaining integrals are now quite easy to deal with after applying the estimate of \eqref{eq2.4.9}:
\begin{equation}
\label{2.4.final}
\begin{aligned}
&\left|E\left[\frac{S_n^{\lambda\enn} S_n^{\mu\enn}}{V^\enn}\right]\right|\\
\leq&
3\epsilon\prod_{j=2}^d\left[\int_\T \left|D(\lambda_1^\enn-\theta_1,v_1^\enn)\right|\cdot\left|D(\mu_1^\enn-\theta_1,v_1^\enn)\right|d|\ma_\T|(e^{i\theta})\right]\\
\leq&3\epsilon
\end{aligned}
\end{equation}

Since $\epsilon$ was arbitrary, the lemma is proved.\qed

\subsection{Another Extension}\label{C2S5}
The hypothesis \eqref{eq2.4.1} from Lemma \ref{lemma2.3}, is obviously essential; if $\lambda^\enn$ and $\mu^\enn$ get ``too close'' to each other (i.e., closer than \eqref{eq2.4.1}), then it should be expected that the asymptotic behavior of the covariance will be closer to situation of Theorem \ref{S1T1}.

That being said, it is only natural to wonder whether Lemma \ref{lemma2.3} couldn't be extended to more than two sequences, so long as the sequences didn't get ``too close''.  Fundamentally, this is the substance of the following result:

\begin{lemma}\label{Lemma2.5}
Let $X$ be a CCWS $\varrho\p$-mixing random field.  Suppose $\{v^\enn\}$ is a sequence that satisfies \eqref{2.1.min}. Let $\{\lambda^\jenn\}_{n=1}^\infty$, $j=1,2,\dots,m$ be a collection of $m$ sequences which converge to $\lambda\in\pid$.  Suppose further that for every pair $j\neq k$, there are a corresponding $\delta(j,k)$, $0<\delta(j,k)<1/2$ and an $N(j,k)\in\N$ such that for every $n\geq N(j,k)$, there is at least one subscript $s$ so that
\begin{equation}
\label{2.5.min}
\left|\lambda_s^{(j,n)}-\lambda_s^{(k,n)}\right|>\frac1{(v_s^\enn)^{1/2-\delta(j,k)}}.
\end{equation}
\begin{enumerate}[\upshape (i)]
\item Then
\[\lim_{n\to\infty}\max_{j\neq k}
\left\{\left|E\left[\frac{S_n^{\lambda(j,n)}\overline{S_n^{\lambda(k,n)}}}{V^\enn}\right]\right|\right\}=0\]
\item If $w\in\Z^d$ is a fixed vector, define $\bwn:=\{k+w:k\in\bn\}$.  With this notation, it holds that
\begin{equation}
\lim_{n\to\infty}\sup_{w\in\Z^d}\max_{j\neq k}\left\{\left|E\left[\frac{S_\bwn^{\lambda\jenn}\overline{S_\bwn^{\lambda(k,\enn)}}}{V^\enn}\right]\right|\right\}=0
\end{equation}
\end{enumerate}
\end{lemma}

\pf To prove (i), first choose $\epsilon>0$.  Let $\delta^*:=\min_{j\neq k}\{\delta(j,k)\}$.  Let $f$ denote the spectral density of $X$ and suppose that $f\leq M$ on $\T^d$.  For any $n$ and any $s$, $1\leq s\leq d$,
\[
\frac1{(v_s^\enn)^{\delta(j,k)}}\leq \frac1{(v_s^\enn)^{\delta^*}}.
\]

Since there are a finite number of pairs $j\neq k$, it is possible to find an $N^*\in\N$ so that \eqref{2.5.min} holds for every pair $j\neq k$ and all $n\geq N^*$.  Since \eqref{2.1.min} holds, we can assume that $N^*$ is also such that for every $n\geq N^*$ and every $s$, $1\leq s\leq d$,
\[
\frac{3M\sqrt\pi}{(v_s^\enn)^{\delta^*}}<\epsilon\]

To complete the proof, pick a pair of indices $j\neq k$, and follow the same steps as in the proof of Lemma \ref{lemma2.3}, where, starting at \eqref{eq248}, replace every $\delta$ with a $\delta^*$.  The bound obtained in this manner is independent of the pair $(j,k)$, hence (i) is proved.

To prove (ii), notice again that if $w\in\Z^d$, then
\[
\begin{aligned}
E\left[X_{k+w}^\lam \overline{X_{j+w}^\mam}\right]&=e^{-iw\cdot\lambda}e^{iw\cdot\mu}E\left[e^{-ik\cdot\lam}X_{k+w} e^{ij\cdot\mu}\overline{X_{j+w}^\mam}\right]\\
=&e^{-iw\cdot\lambda}e^{iw\cdot\mu}E\left[X_k^\lam \overline{X_j^\mam}\right]
\end{aligned}
\]
where the last equality follows from the stationarity properties.  Thus
\[
\left|E\left[S_\bwn^\lam S_\bwn^\mam\right]\right|=\left|E\left[S_n^\lam S_n^\mam\right]\right|,
\]
and the result follows.\qed

\subsection{``Miller's'' Lemma}\label{C2S6}
The result of this section is named after Curtis Miller, who proved a similar result in \cite{6} under different conditions.  The result will be useful when using the Cramer-Wold device (Theorem \ref{Cramer}) in the next section.

A bit of notation is necessary before the theorem can be stated.  Let $\bv:=(a_1,b_1,\dots,a_n,b_n)\in\R^{2m}$ and let $\textbf{z}:=(z_1,z_2,\dots,z_m)\in\C^m$.  Define the function
\begin{equation}
\label{funcG}
G(\bv,\textbf z)=\sum_{j=1}^m\left(a_j\Re z_j+b_j\Im z_j\right)
\end{equation}
Also, let $\|\bv\|^2=a_1^2+b_1^2+\dots+a_m^2+b_m^2$

As before, if $w\in\Z^d$, let $\bwn:=\{k+w:k\in\bn\}$.  Then, if $\{\lambda^\jenn\}_{n=1}^\infty$, $j=1,2,\dots,m$ are all sequences of elements of \pid, define the vector
\begin{equation}
\label{vecS}
\bS_n^w\left(\lambda^{(1,n)},\lambda^{(2,n)}\dots,\lambda^{(m,n)}\right):=
\left(S_\bwn^{\lambda(1,n)},S_\bwn^{\lambda(2,n)},\dots,S_\bwn^{\lambda(m,n)}\right)
\end{equation}
In the context of Lemma \ref{LemmaM} below, the sequences $\{\lambda^\jenn\}$ will be chosen to satisfy certain conditions, but will otherwise remain fixed.  As usual, the sequence $\{v^\enn\}$ will be assumed to satisfy \eqref{2.1.min}, but it will also be fixed.  Thus, a shorthand version of \eqref{vecS} will be possible:
\[
\bS_n^w:=\bS_n^w\left(\lambda^{(1,n)},\lambda^{(2,n)}\dots,\lambda^{(m,n)}\right).
\]
Also, in the context below, the vector $\bv$ will be assumed fixed, and so instead of writing $G\left(\bv,\bS_n^w\right)$ we shall simply write $G\left(\bS_n^w\right)$
\begin{lemma}[``Miller's Lemma'']\label{LemmaM}
Let $X$ be a CCSS random field which is $\varrho\p$-mixing.  Assume $E|X_0|^2=\sigma^2<\infty$.  Let $f(e^{i\lambda})=f\lam$ denote the (continuous) spectral density of $X$.  Suppose that $\lambda\in\at$ and let $\{\lambda^\jenn\}_{n=1}^\infty$, $j=1,2,\dots,m$ be sequences which converge to $\lambda$ and also satisfy the conditions of Lemma \ref{Lemma2.5}.  Suppose $\{v^\enn\}$ is a sequence which satisfies \eqref{2.1.min}.  Let $\emph{\bv}=(a_1,b_1,\dots,a_n,b_n)\in\R^{2m}$ be an arbitrary (but fixed) vector.  Then
\begin{equation}
\label{M1}
\lim_{n\to\infty}\sup_{w\in\Z^d}
\left|
\frac12f\lam\|\emph{\bv}\|^2-\frac1{V^\enn}E\left[G\left(\bS_n^w\right)\right]^2
\right|=0.
\end{equation}
\end{lemma}

\pf If $\bv=0$, there is nothing to show.  So, assume the contrary.

Let $1>\epsilon>0$.  Since $f$ is continuous on the compact set
$[-\pi,\pi]^d$, it is uniformly continuous.  Therefore, there exists a
$\delta>0$ such that if $\mu$ satisfies $|\lambda-\mu|<\delta$, then

\begin{equation}
\label{Miller1}
|f\lam-f(\mu)|\leq\frac{\epsilon}{8m\|{\bv}\|^2}.
\end{equation}

Next, Theorem \ref{S1T1} implies that
\[
\lim_{n\to\infty}EI_n^\mu=f(\mu)
\]
uniformly over all
$\mu\in[-\pi,\pi]^d.$  There is therefore an $N_1$ such
that for all $n>N_1$ and all $\mu\in[-\pi,\pi]^d$,
\begin{equation}
\label{Miller1.1}
\left|EI_n^\mu-f(\mu)\right|\leq\frac{\epsilon}{8m\|\bv\|^2}.
\end{equation}

Next, consider the compact set
$[-\pi,\pi]^d\supset K:=\left(\bigcup_{j,n}\lambda^{(j,n)}\right)\cup\{\lambda\}$. Since
$\lambda\in\mathfrak{P}$, $(\lambda_s)^2\neq1$
for some subscript $s$, $1\leq s\leq d$.  Since there are only a
finite number of sequences $\{\lambda^{(j,n)}\}$, $j=1,2,\dots,m$, there
exists an $N_2>0$ such that for all $n>N_2$,
$\lambda_s^{(j,n)}\cdot\lambda_s^{(k,n)}\neq1$ for all pairs $j$ and
$k$ (note that this holds even when $j=k$).  Eliminating a finite collection
of points from $K$ (that is, all the $\lambda^{(j,n)}$ where $n\leq
N_2$), leaves another compact set. So, apply Lemma
\ref{lemma2.2} to find an $N_3\geq N_2$ such that
\begin{equation}
\label{Miller2}
\text{for all }n>N_3,\smsp
\left|E\left[\frac{S_n^{\lambda(j,n)}S_n^{\lambda(k,n)}}{V^\enn}\right]\right|<\frac{\epsilon}{8m\|\bv\|^2}.
\end{equation}
It now follows that
\begin{equation}
\label{Miller3} \text{for all }n>\max\{N_1,N_3\}, \smsp
\left|E\left[\frac{S_n^{\lambda(j,n)}S_n^{\lambda(k,n)}}{V^\enn}\right]\right|<\frac{\epsilon}{8m\|\bv\|^2}
\end{equation}
which holds for all pairs $j$ and $k$, including when $j=k$

Lemma \ref{Lemma2.5} also applies to the current context.  Therefore, there
exists an $N_4$ such that
\begin{equation}
\label{Miller4}
\text{for all }n>N_4,\smsp
\left|E\left[\frac{S_n^{\lambda(j,n)}\overline{S_n^{\lambda(k,n)}}}{V^\enn}\right]\right|\leq\frac{\epsilon}{8m\|\bv\|^2},
\end{equation}
which holds for every pair $j\neq k$.
Claim \ref{Claim0} and equation \eqref{eq2.3.b} imply that the bounds \eqref{Miller3} and \eqref{Miller4} will also hold for the $S_\bwn^\lam$'s as well.

Now let
\[
\begin{aligned}
R_n^{(\mu)}&:=\Re S_n^{(\mu)}\\
Q_n^{(\mu)}&:=\Im S_n^{(\mu)}.
\end{aligned}
\]
Then:
\begin{equation}
\begin{aligned}
\label{Miller5}
\Re\left\{E\frac{S_n^{\lambda(j,n)}S_n^{\lambda(k,n)}}{V^\enn}\right\}
    =E\frac{R_n^{\lambda(j,n)}R_n^{\lambda(k,n)}}{V^\enn}-E\frac{Q_n^{\lambda(j,n)}Q_n^{\lambda(k,n)}}{V^\enn}
\end{aligned}
\end{equation}
\begin{equation}
\begin{aligned}
\label{Miller5.1}
\Im\left\{E\frac{S_n^{\lambda(j,n)}S_n^{\lambda(k,n)}}{V^\enn}\right\}=
E\frac{R_n^{\lambda(j,n)}Q_n^{\lambda(k,n)}}{V^\enn}+E\frac{Q_n^{\lambda(j,n)}R_n^{\lambda(k,n)}}{V^\enn}
\end{aligned}
\end{equation}

\begin{equation}
\begin{aligned}
\label{Miller5.2}
\Re\left\{E\frac{S_n^{(\lambda(j,n))}\overline{S_n^{\lambda(k,n)}}}{V^\enn}\right\}=
E\frac{R_n^{\lambda(j,n)}R_n^{\lambda(k,n)}}{V^\enn}+E\frac{Q_n^{\lambda(j,n)}Q_n^{\lambda(k,n)}}{V^\enn}
\end{aligned}
\end{equation}
\begin{equation}
\begin{aligned}
\label{Miller5.3}
\Im\left\{E\frac{S_n^{(\lambda(j,n))}\overline{S_n^{\lambda(k,n)}}}{V^\enn}\right\}
=E\frac{Q_n^{\lambda(j,n)}R_n^{\lambda(k,n)}}{V^\enn}-E\frac{R_n^{\lambda(j,n)}Q_n^{\lambda(k,n)}}{V^\enn}
\end{aligned}
\end{equation}
The equations \eqref{Miller3} and \eqref{Miller4} imply that, when $j\neq k$ and
$n\geq\max\{N_1,N_3,N_4\}$, the modulus of each of the terms
\eqref{Miller5}-\eqref{Miller5.3} must be less than $\epsilon/8m\|\bv\|^2$. The following claim allows this bound to be applied to the four terms on the right-hand sides:
\begin{claim}\label{Claim1} If $x,y\in\C$, $|x+y|\leq c$ and $|x-y|\leq c$, then
$|x|< c$ and $|y|< c$.\end{claim}
\pf Heuristically, at least one of the diagonals in a parallelogram is longer than every side. \qed\\

Claim \ref{Claim1}, equations \eqref{Miller5}-\eqref{Miller5.3}, and the remarks just following \eqref{Miller5}-\eqref{Miller5.3}, imply that, when $j\neq k$:
\begin{equation}
\label{Miller5.4}
\left|E\frac{R_n^{\lambda(j,n)}R_n^{\lambda(k,n)}}{V^\enn}\right|\leq\frac{\epsilon}{8m\|\bv\|^2}
\end{equation}
with the same bound holding for the three other terms like it from right-hand sides of \eqref{Miller5}-\eqref{Miller5.3}.  Because the bounds in \eqref{Miller3} and \eqref{Miller4} are uniform over all $j\neq k$, so is the bound in \eqref{Miller5.4} (and implicitly the bounds on the other three terms).

When $j=k$, the situation is a bit different:
\[
\begin{aligned}
E\frac{\left(S_n^{\lambda(j,n)}\right)^2}{V^\enn}&=\underbrace{E\frac{\left(R_n^{\lambda(j,n)}\right)^2}{V^\enn}-E\frac{\left(Q_n^{\lambda(j,n)}\right)^2}{V^\enn}}_{(i)}
    +\underbrace{2iE\frac{R_n^{\lambda(j,n)}Q_n^{\lambda(j,n)}}{V^\enn}}_{(ii)}\\
E\frac{\left|S_n^{\lambda(j,n)}\right|^2}{V^\enn}&=\underbrace{E\frac{\left(R_n^{\lambda(j,n)}\right)^2}{V^\enn}+E\frac{\left(Q_n^{\lambda(j,n)}\right)^2}{V^\enn}}_{(iii)}
\end{aligned}
\]
Lemma \ref{lemma2.2} and equation \eqref{Miller3} guarantee that terms (i) and (ii) above each
converge to zero uniformly over all $j$.
Term (iii), however, is (by \eqref{Miller1.1}), close to
$f(\lambda^{(j,n)})$.  These two facts together imply that
\begin{equation}
\begin{aligned}
\label{Miller6}
\left|E\frac{\left(R_n^{\lambda(j,n)}\right)^2}{V^\enn}-\frac{f(\lambda^{(j,n)})}2\right|&<\frac\epsilon{8m\|\bv\|^2} \quad\text{ and }\\
\left|E\frac{\left(Q_n^{\lambda(j,n)}\right)^2}{V^\enn}-\frac{f(\lambda^{(j,n)})}2\right|&<\frac\epsilon{8m\|\bv\|^2}.
\end{aligned}
\end{equation}
Again, these bounds hold uniformly over all $j$.

We now have need of the following claim:

\begin{claim}\label{Claim2}
If $x_j$, $j=1,2,\dots,n$ are real numbers, then $\left(\sum_{j=1}^n|x_j|\right)^2\leq n \sum_{j=1}^nx_j^2$.
\end{claim}\pf
\[
\begin{aligned}
\left(\sum_{j=1}^n|x_j|\right)^2&=\sum_{j=1}^nx_j^2+\sum_{j\neq k}|x_jx_k|\\
&\leq\sum_{j=1}^nx_j^2+\sum_{j\neq k}\frac{x_j^2+x_k^2}2=\sum_{j=1}^nx_j^2+(n-1)\sum_{j=1}^nx_j^2\\
&=n\sum_{j=1}^nx_j^2.
\end{aligned}
\]\qed

Now
\begin{equation}
\begin{aligned}
\label{Miller7}
\frac{E\left[G\left(\bS_n^0\right)\right]^2}{V^\enn}
=\frac1{V^\enn}\sum_{j=1}^m\left(a_j^2E\left[\left(R_n^{\lambda(j,n)}\right)^2\right]+b_j^2E\left[\left(Q_n^{\lambda(j,n)}\right)^2\right]\right)
+g_n^\epsilon
\end{aligned}
\end{equation}
where $g_n^\epsilon$ is a ``garbage'' term that includes all the
missing terms.  The proof of Claim \ref{Claim2}, together with equation \eqref{Miller5.4} prove that

\begin{equation}
\label{Miller8.1}
\begin{aligned}
|g_n^\epsilon|&\leq \sum_{j\neq k}
\left(|a_jb_k|\frac\epsilon{8m\|\bv\|^2}\right)\\
&\leq \epsilon\frac{(2m-1)}{8m\|\bv\|^2}\cdot\|\bv\|^2\\
&< \frac\epsilon4
\end{aligned}
\end{equation}

Hence, if the inequalities of \eqref{Miller8.1}, \eqref{Miller6} and \eqref{Miller1.1} hold, then
\begin{equation}
\begin{aligned}
\label{Miller9}
&\left|\frac12f\lam\sum_{j=1}^m(a_j^2+b_j^2)-\frac1{V^\enn}E\left[G\left(\overrightarrow{S}_n^\nu\right)\right]^2\right|\\
\leq&\smsp|g_n^\epsilon|+\left|\frac12\sum_{j=1}^m(a_j^2+b_j^2)\left(f\lam-f(\lambda^{(j,n)})\right)\right|\\
    &+\left|\sum_{j=1}^ma_j^2\left(\frac{f(\lambda^{(j,n)})}{2}-\frac{E\left[\left(R_n^{\lambda(j,n)}\right)^2\right]}{V^\enn}\right)\right|\\
    &+\left|\sum_{j=1}^mb_j^2\left(\frac{f(\lambda^{(j,n)})}{2}-\frac{E\left[\left(Q_n^{\lambda(j,n)}\right)^2\right]}{V^\enn}\right)\right|\\
\leq&\epsilon,
\end{aligned}
\end{equation}
which proves the lemma. \qed

\newpage
\section{The Main Central Limit Theorem}\label{C3}

\subsection{Statement of the Main Result}\label{C3S1}

The goal of this section will be to prove the following Central Limit Theorem.  (For the definition of $\bS_n^0$, see \S\ref{C2S6}; for the definition of $\at$, see \S\ref{C2S1})

\begin{theorem}\label{N3}
Let $d$ be a positive integer and suppose $X:=\{X_k,k\in\Z^d\}$ is a $\varrho\p$-mixing, CCSS random field such that $E|X_k|^2=\sigma^2<\infty$.  Let $f\lam:=f(e^{i\lambda})$ be the (continuous) spectral density of $X$. Let $\lambda\in\at$, and let $\{\lambda^\jenn\}_{n=1}^\infty$, $j=1,2,\dots,m$ be sequences of elements of \pid which converge to $\lambda$, and which satisfy the conditions of Lemma \ref{Lemma2.5}. Suppose $\{v^\enn\}$ is a sequence of vectors that satisfies \eqref{2.1.min}, i.e.
\[\lim_{n\to\infty}\min\{v_1^\enn,v_2^\enn,\dots,v_d^\enn\}=\infty.
\] Then
\begin{equation}\label{Main3}
    \frac{\bS_n^0}{\sqrt{V^\enn}}\Rightarrow \textbf{Z}
\end{equation}
where $\textbf{Z}:\Omega\rightarrow \R^{2m}$ has the normal distribution with the $2m\times2m$ covariance matrix
\[\Upsilon_m^\lam:=
\begin{bmatrix}
f\lam&0&\dots&0\\
0&f\lam&\dots&0\\
\vdots&\vdots&\ddots&\vdots\\
0&0&\dots&f\lam
\end{bmatrix}.
\]
\end{theorem}

Generally speaking, the proof consists of two reductions.  The first involves truncating the individual random variables (i.e., the $X_k$'s).  The second reduction involves the Bernstien blocking argument.  There, a number of conveniently-selected truncated random variables will be eliminated from the normed sums.  Since the total number of these removed random variables will be small when compared to the total number of summands they will be shown to be negligible insofar as weak convergence is concerned.

In both reductions, the following lemma will play a significant role:

\begin{lemma}\label{Lemma3.1}
Let $\{Y_n\}_{n=1}^\infty$ be a sequence of centered random variables on the same probability space and suppose $E|Y_n|^2\to0$ as $n\to\infty$.  Then $Y_n\Rightarrow0$.
\end{lemma}

\pf Choose $\epsilon>0$. Apply Chebychev's inequality:
\[P\left(|Y_n-0|>\epsilon\right)\leq\frac1{\epsilon^2}Var(Y_n)\leq\frac1{\epsilon^2}E|Y_n|^2\to0.\] But convergence in probability implies weak convergence (see \cite{1}, Theorem 25.2).  \qed

This lemma, when combined with Theorem \ref{slutsky}, will provide the means to show that both reductions are valid.\\

\noindent\textbf{Start of Proof of Theorem \ref{N3}:}  Fix $\lambda\in\at$.  If $f\lam=0$, apply the Cramer-Wold Device (Theorem \ref{Cramer}). Then, Lemmas \ref{LemmaM} and \ref{Lemma3.1} imply that $\bS_n^0\Rightarrow0$.  Therefore, assume $f\lam>0$.\\

The proof of Theorem \ref{N3} is officially concluded in \S\ref{C3S7}.  The intervening subsections should be considered as parts of the proof, which has been subdivided for clarity.

\subsection{The Truncated Random Variables}\label{C3S2}
It might prove useful to the reader to review the notations of \S\ref{C2S6}.

Here and below, whenever $k\in\Z^d$, denote by $\bk=\prod_{i=1}^d k_i$.  Then, let $0<q$ be a real number, and suppose that $\mu\in\pid$.  For each $k\in\Z^d$, define

\[
\begin{aligned}
B_{k,q}^{(\mu)}:&=X_k^{(\mu)}\mathbb{I}_{\{|X_k|\leq\bk^q\}}-E\left[X_k^{(\mu)}\mathbb{I}_{\{|X_k|\leq\bk^q\}}\right]\\
T_{k,q}^{(\mu)}:&=X_k^{(\mu)}\mathbb{I}_{\{|X_k|>\bk^q\}}-E\left[X_k^{(\mu)}\mathbb{I}_{\{|X_k|>\bk^q\}}\right].
\end{aligned}
\]
Notice that
$E[B_{k,q}^{(\mu)}]=E[T_{k,q}^{(\mu)}]=0$ for all $k$, $q$, and
$\mu$.  Furthermore, for all $\mu,\nu\in\pid$, and every $\omega\in\Omega$
\begin{equation}\label{eq3.2.1}
    \left|B_{k,q}^{(\mu)}(\omega)\right|=\left|B_{k,q}^{(\nu)}(\omega)\right|
    \quad\&\quad\left|T_{k,q}^{(\mu)}(\omega)\right|=\left|T_{k,q}^{(\nu)}(\omega)\right|.
\end{equation}
The equations in \eqref{eq3.2.1} hold because complex constants can be factored out of the expectation.

Next, define
\[
\begin{aligned}
S_{n,q}^{(\mu)}:&=\sum_{k\in\bn}B_{k,q}^{(\mu)}\\
Z_{n,q}^{(\mu)}:&=\sum_{k\in\bn}T_{k,q}^{(\mu)}\\
R_{n,q}^{(\mu)}:&=\Re[S_{n,q}^{(\mu)}]\\
Q_{n,q}^{(\mu)}:&=\Im[S_{n,q}^{(\mu)}]\\\\
\end{aligned}
\]

As before, if $w\in\Z^d$ is fixed, define $\bwn:=\{k+w:k\in\bn\}$.  Then, let

\[
S_{n,q}^{\mu,w}:=\sum_{k\in\bwn}B_{k,q}^\mam.
\]
Now let
\begin{equation}\label{eq3.1.1}
    \bS_{n,q}^w\Big(\lambda^{(1,n)},\lambda^{(2,n)},\dots,\lambda^{(m,n)}\Big)
    :=\left(S_{n,q}^{\lambda(1,n),w},S_{n,q}^{\lambda(2,n),w},\dots,S_{n,q}^{\lambda(m,n),w}\right)
\end{equation}
which is a $\C^m$-valued random vector.  The vectors $\lambda^\jenn$ in the current context are assumed to satisfy certain conditions, but otherwise they are assumed to be fixed.  Therefore, the following shorthand notation for \eqref{eq3.1.1} is possible:
\[\bS_{n,q}^w:=\bS_{n,q}^w\Big(\lambda^{(1,n)},\lambda^{(2,n)},\dots,\lambda^{(m,n)}\Big).\]

\subsection{Reduction I: A CLT for the Truncated Sums}\label{C3S3}
The proceeding lemma is the first reduction mentioned in \S\ref{C3S1}:

\begin{lemma}\label{qLemma}
In the same context and with the same notations as Theorem \ref{N3}, and for any $q$, $0<q<\frac14$.
\begin{equation}\label{Main3.2}
    \frac{\left(R_{n,q}^{\lambda(1,n)},Q_{n,q}^{\lambda(1,n)},\dots,R_{n,q}^{\lambda(m,n)},Q_{n,q}^{\lambda(m,n)}\right)}{\sqrt{V^\enn}}\Rightarrow N\left(0,\Upsilon_m^\lam\right).
\end{equation}
Additionally,
\begin{equation}\label{Secondary3.2}
\frac1{\sqrt{V^\enn}}\left(Z_{n,q}^{\lambda(1,n)},Z_{n,q}^{\lambda(2,n)},\dots,Z_{n,q}^{\lambda(m,n)}\right)\Rightarrow0\in\C^m.\vspace{.2in}
\end{equation}
\end{lemma}

\noindent The proof of \eqref{Main3.2} is given in \S\ref{C3S7}.\\

\noindent \textbf{Proof of Lemma \ref{qLemma}, Equation \eqref{Secondary3.2}:} Let $a=(a_1,a_2,\dots,a_m)\in\R^m$
be arbitrary but fixed. To prove \eqref{Secondary3.2}, it is enough to
show (see \cite{1.5}, p.27) the
following two statements:
\begin{equation}
\label{qLem1}
\frac1{\sqrt{\Lambda_n}}\sum_{j=1}^ma_j\Re\left[Z_{k,q}^{\lambda(j,n)}\right]\Rightarrow
0\in\R;
\end{equation}

\begin{equation}
\label{qLem2}
\frac1{\sqrt{\Lambda_n}}\sum_{j=1}^ma_j\Im\left[Z_{k,q}^{\lambda(j,n)}\right]\Rightarrow
0\in\R.
\end{equation}
This follows from the Cramer-Wold Device, and the fact that if $X_n$
and $Y_n$ are $\R^m$-valued random vectors which both converge
weakly to zero in $\R^m$, then $X_n+iY_n$ converges weakly to zero
(considered as an element of $\C^m$).  The proofs of \eqref{qLem1}
and \eqref{qLem2} are analogous, so only the proof of
\eqref{qLem1} is presented below.

Define the random field $Y_n^{(a)}:=\left\{Y_{k,n}^{(a)}=\sum_{j=1}^m a_j\Re T_{k,q}^{\lambda\jenn}:k\in\Z^d\right\}$. Since $X$ is $\varrho\p$-mixing, so is $Y_n^{(a)}$.  Moreover, $\varrho(Y_n^{(a)},n)\leq \varrho(X,n)$. Therefore, the Rosenthal Inequality (Theorem \ref{Brad2}) with $\beta=2$ implies that
\begin{equation}\label{eq3.3.5}
E\left(\sum_{k\in\bn}\frac{Y_{k,n}^{(a)}}{\sqrt{V^\enn}}\right)^2\leq\frac C{V^\enn}\sum_{k\in\bn} E\left(Y_{k,n}^{(a)}\right)^2.
\end{equation}
Here, the constant $C$ is the one associated with the original field $X$, not a new one which might depend on the field $Y_n^{(a)}$.
Notice that Claim \ref{Claim2} implies that
\begin{equation}\label{eq3.3.6}
E\left(Y_{k,n}^{(a)}\right)^2\leq m\sum_{j=1}^ma_j^2E\left(\Re\left[T_{k,q}^{\lambda\jenn}\right]\right)^2
\leq m\sum_{j=1}^ma_j^2E\left|T_{k,q}^{\lambda\jenn}\right|^2.
\end{equation}
Equation \eqref{eq3.2.1} implies that
\begin{equation}\label{eq3.3.7}
    \sum_{j=1}^ma_j^2E\left|T_{k,q}^{\lambda\jenn}\right|^2= m\|a\|^2E\left|T_{k,q}^{\lambda(1,n)}\right|^2.
\end{equation}
Now put \eqref{eq3.3.5}, \eqref{eq3.3.6}, and \eqref{eq3.3.7} together:
\[
E\left(\frac1{\sqrt{\Lambda_n}}\sum_{j=1}^ma_j\Re\left[Z_{k,q}^{\lambda(j,n)}\right]\right)^2\leq
\frac {Cm^2\|a\|^2}{V^\enn}\sum_{k\in\bn} E\left|T_{k,q}^{\lambda(1,n)}\right|^2.
\]
Since $E|X_k|^2=\sigma^2<\infty$ for all $k$, the values $\left\{E\left|T_{k,q}^{\lambda(1,n)}\right|^2\right\}$, $k\in\Z^d$, satisfy the conditions of Lemma \ref{LemmaZ}.  Hence
\[
E\left(\frac1{\sqrt{\Lambda_n}}\sum_{j=1}^ma_j\Re\left[Z_{k,q}^{\lambda(j,n)}\right]\right)^2\to0.
\]
Now apply Lemma \ref{Lemma3.1} to finish the proof of \eqref{Secondary3.2}.

The proof of Lemma \ref{qLemma} continues into the next sections, concluding in \S\ref{C3S7}.

\subsection{The Bernstein Blocking Technique}\label{C3S4}
What follows is an argument that involves the Bernstein Blocking technique (which actually dates back to at least Markov).  Heuristically, the gist of the argument involves slicing $\bS_{n,q}^0$ like a loaf of bread, except that the ``width'' of the slices is not to be uniform.  Instead, the first slice should be thick, the second thin, the third thick, the fourth thin, and so on.  The thin slices should grow in thickness as $\bS_{n,q}^0$ grows in size, and since they come between the thick slices, these latter pieces should be quasi-independent because of the mixing condition.  However, when taken all together, the thin slices can't account for too much of $\bS_{n,q}^0$, since when $\bS_{n,q}^0$ is normed, it will be desirable to be able to neglect the small slices and focus on the big ones.  Thus, much care must be taken to ensure that both of these criteria are satisfied.\\

Many of the notations from \S\S\ref{C2S6} and \ref{C3S2} will be used throughout this section.

Let $\bv=(a_1,b_1,a_2,b_2,\dots,a_m,b_m)\in\R^{2m}$ be arbitrary but fixed and denote by $\|\bv\|$ the Euclidean norm of $\bv$. Let $G(\bv,\textbf{z})$ be as in \eqref{funcG}.  The Cramer-Wold device implies that if it can be shown that
\begin{equation}\label{eq3.4.1}
    \frac{G\left(\bv,S_{n,q}^0\right)}{\sqrt{V^\enn}}\Rightarrow N\left(0,f\lam\|\bv\|^2\right)
\end{equation}
then Lemma \ref{qLemma} follows.  (It is a standard fact of probability theory that a linear combination of jointly normal random variables is again normal.  Hence the right-hand side of \eqref{eq3.4.1} follows simply by calculating the variance of that particular linear combination.)

To that end, let $\lfloor\cdot\rfloor$ denote the largest integer less than or equal to the argument, and define
\begin{equation}\label{eq3.4.2}
    s\enn:=\left\lfloor \sqrt[3]{v_1^\enn}\right\rfloor
\end{equation}
($s\enn$ will be the width of the small ``slices'' or ``blocks'').  Since $v^\enn$ satisfies \eqref{2.1.min}, $s\enn\to\infty$ as $n\to\infty$.
Now define
\begin{equation}\label{eq3.4.3}
    p\enn:=\min\left\{s\enn,\left\lfloor\frac1{\sqrt{\rho\p(s\enn)}}\right\rfloor\right\}
\end{equation}
($p\enn$ will determine the number of slices made).  Again, notice that $p\enn\to\infty$ with $n$.  Finally, define $r\enn$ to be the positive integer which satisfies\\
\begin{equation}\label{eq3.4.4}
    (r\enn-1+s\enn)p\enn\leq v_1^\enn<(r\enn+s\enn)p\enn\\\\
\end{equation}
($r\enn$ is the width of the big ``slices'').  Once again, it follows easily from \eqref{eq3.4.2} and \eqref{eq3.4.3} that $r\enn\to\infty$ with $n$.

There are a couple of important inequalities that fall out directly from \eqref{eq3.4.2}-\eqref{eq3.4.4}.  First of all,
\[
v_1^\enn-\left(v_1^\enn\right)^{2/3}\leq v_1^\enn-p\enn s\enn\leq p\enn r\enn,
\]
and therefore:
\begin{equation}\label{eq3.4.0}
    v_1^\enn-p\enn r\enn\leq \left(v_1^\enn\right)^{2/3},
\end{equation} and so a simple argument yields
\begin{equation}\label{eq3.4.0b}
    \lim_{n\to\infty}\frac{p\enn r\enn}{v_1^\enn}=1.
\end{equation}

It's time to construct the thick slices.  For every integer $l$, $1\leq l\leq p\enn$, and every $n$, define the subset of $\Z^d$

\begin{equation}\label{boxes}
\mathbb{B}(l,n):=\left\{k\in\bn:(l-1)(r\enn+s\enn)< k_1\leq lr\enn +(l-1)s\enn\right\}\vspace{.2in}
\end{equation}

Notice that $\Card\left(\mathbb{B}(l,n)\right)=r\enn\times v_2^\enn\times\cdots \times v_d^\enn$.
Also define the random variables
\begin{equation}\label{Gamma}
\GG=\Gamma(l,n,q):=\frac{\sum_{k\in\mathbb{B}(l,n)}
G\left(\bv,\left(B_{k,q}^{\lambda(1,n)},\dots,B_{k,q}^{\lambda(m,n)}\right)\right)}{\sqrt{V^\enn}}.
\end{equation}
The \GG's are the ``thick'' slices.  It will also be helpful to collect the tails from these thick slices.  To do so, define for each $l$, $1\leq l\leq p\enn$, the random variable

\begin{equation}\label{xis}
    \xi(l,n)=\xi(l,n,q):=\frac{\sum_{k\in\mathbb{B}(l,n)}
G\left(\bv,\left(T_{k,q}^{\lambda(1,n)},\dots,T_{k,q}^{\lambda(m,n)}\right)\right)}{\sqrt{V^\enn}}.
\end{equation}
We will come back to the $\Gamma$'s and $\xi$'s in a moment.

As mentioned earlier, it is now important to show that the small slices are negligible. Instead of dealing with each small slice individually, it is possible to take them all together; consider the set
\[
\ZZ\enn:=\bn\setminus\left(\mathbb{B}(1,n)\cup\mathbb{B}(2,n)\cup\dots\cup\mathbb{B}(m,n)\right).
\]
The cardinality of $\ZZ\enn$ is equal to (\emph{cf.} equation \eqref{eq3.4.0})

\begin{equation}\label{cardz}
\left(v_1^\enn-p\enn r\enn\right)\cdot v_2^\enn\cdots v_d^\enn\leq\left(v_1^\enn\right)^{2/3}\cdot v_2^\enn\cdots v_d^\enn=o(V^\enn)
\end{equation}
It is now desirable to show the following

\begin{lemma}\label{Lemma3.4}  With all the notations and assumptions of the current context,
\begin{equation}\label{main3.4}
\frac{\sum_{k\in\ZZ\enn}G\left(\bv,\left(B_{k,q}^{\lambda(1,n)},\dots,B_{k,q}^{\lambda(m,n)}\right)\right)}{\sqrt{V^\enn}}
\Rightarrow0\in\R.
\end{equation}

\end{lemma}

\noindent\textbf{Proof of Lemma \ref{Lemma3.4}:}  Suppose the following statement can be proved:
\begin{equation}\label{second3.4}
    \frac{\sum_{k\in\ZZ\enn}
    a_1\Re B_{k,n}^{\lambda(1,n)}+
    b_1\Im B_{k,n}^{\lambda(1,n)}}{\sqrt{V^\enn}}\Rightarrow0.
\end{equation}
Now, there is nothing particularly interesting about the subscripts ``1'' in \eqref{second3.4}; the proof below does not rely on the the subscript.  That is to say, if \eqref{second3.4} holds (i.e., when $j=1$), then it also holds for every other $j$, $2\leq j\leq m$.  This, combined with Slutsky's Theorem (Theorem \ref{slutsky}), shows that \eqref{main3.4} holds as well.

So, to prove \eqref{second3.4}, first assume that at least one of $a_1$ or $b_1$ is non-zero, since otherwise there is nothing to show.
Now consider the variance:
\[
M_2^{\enn}:=E\left[\left(\frac{\sum_{k\in\ZZ\enn}
    a_1\Re B_{k,n}^{\lambda(1,n)}+
    b_1\Im B_{k,n}^{\lambda(1,n)}}{\sqrt{V^\enn}}\right)^2\right].
\]
The random field $\left\{a_1\Re B_{k,n}^{\lambda(1,n)}+b_1\Im B_{k,n}^{\lambda(1,n)}:k\in\Z^d\right\}$ is $\rho\p$-mixing, hence by the Rosenthal inequality
\begin{equation}\label{eq3.4.8}
M_2^{\enn}\leq C\cdot\frac{\sum_{k\in\ZZ\enn}E\left(
    a_1\Re B_{k,n}^{\lambda(1,n)}+
    b_1\Im B_{k,n}^{\lambda(1,n)}\right)^2}{{V^\enn}}\smsp,
\end{equation}
where again, $C$ is a constant associated with the original field $X$.  Claim \ref{Claim2} implies that

\begin{equation}\label{eq3.4.10}
\begin{aligned}
E\left(
    a_1\Re B_{k,n}^{\lambda(1,n)}+
    b_1\Im B_{k,n}^{\lambda(1,n)}\right)^2&\leq2\max\{a_1^2,b_1^2\}\cdot E\left|B_{k,n}^{\lambda(1,n)}\right|^2\\
    &\leq2 \max\{a_1^2,b_1^2\}\cdot E\left|X_0\right|^2.
\end{aligned}
\end{equation}
Combining \eqref{eq3.4.8} and \eqref{eq3.4.10},
\begin{equation}\label{eq3.4.11}
    M_2^\enn\leq 2C\max\{a_1^2,b_1^2\}E|X_0|^2\frac{\Card(\ZZ\enn)}{V^\enn}\to0,
\end{equation}
where the convergence follows from the cardinality of $\ZZ\enn$ given in \eqref{cardz}.  The statement in  \eqref{second3.4} is therefore true because of Lemma \ref{Lemma3.1}. This completes the proof of Lemma \ref{Lemma3.4}. \qed

\subsection{Reduction II: Independent Blocks}\label{C3S5}
In this section, the goal will be to show that the big ``slices'' of the Bernstien blocking technique  from \S\ref{C3S4} can be thought of as being independent, at least insofar as a Central Limit Theorem is concerned.  The reason that this will be possible is because the big slices are separated by the smaller ones, thus permitting the use of the correlation coefficient to analyze the characteristic function.\\

Notice first of all that the random field  \[\left\{G\left(\bv,\left(B_{k,q}^{\lambda(1,n)},B_{k,q}^{\lambda(2,n)},\dots,B_{k,q}^{\lambda(m,n)}\right)\right):k\in\Z^d\right\}\]
is $\varrho\p$-mixing.  Therefore, if $\mathcal{B}(l,n)$ denotes the $\sigma$-field generated by the $X_k$'s in $\mathbb{B}(l,n)$, i.e.
\[
\mathcal{B}(l,n)=\sigma\left(X_k:k\in\mathbb{B}(l,n)\right),
\]
then, whenever $1\leq l_1\neq l_2\leq p^\enn$,
\begin{equation}\label{eq3.5.0}
\varrho\p\big(\mathcal{B}(l_1,n),\mathcal{B}(l_2,n)\big)\leq \varrho\p(X,s^\enn).
\end{equation}

Denote the characteristic function of the sum of the (dependent) big blocks by $\phi^\enn(t)$:
\[\phi^\enn(t):=E\left[\exp\left\{it\sum_l \GG\right\}\right].\]
Now construct independent copies of the $\GG$'s: let $\DD$ have the same distribution as $\GG$, but suppose that $\Delta(l_1,n)$ and $\Delta(l_2,n)$ are independent whenever $1\leq l_1\neq l_2\leq p^\enn$.  Let their characteristic function be denoted by
\[
\psi^\enn(t):=E\left[\exp\left\{it\sum_l \DD\right\}\right].
\]

\begin{lemma}\label{Lemma3.5}
If $\phi^\enn(t)$ and $\psi^\enn(t)$ are defined as above, then, for all $t\in\R$
\begin{equation}\label{eq3.5.1}
    \lim_{n\to\infty}\left|\phi^\enn(t)-\psi^\enn(t)\right|=0.
\end{equation}
\end{lemma}
Before moving into the proof, notice that a simple application of Theorem \ref{Bill26.3} together with Lemma \ref{Lemma3.5} implies that if $\sum_l\DD$ converges weakly to some distribution, then $\sum_l\GG$ converges weakly to that same distribution.  This will make the second reduction possible, since it means that in order to analyze the end behavior of $\sum_l\GG$, it is sufficient to analyze the end behavior of $\sum_l\DD$.\\

\noindent{\textbf{Proof of Lemma \ref{Lemma3.5}:}}  Equation \eqref{eq3.5.0} and the paper by Withers \cite{12} (which is needed to show that the dependence coefficient can be applied to complex-valued functions functions) imply that the following inequality holds:

\[
\begin{aligned}
    &\left|\phi^\enn(t)-E\big[\exp\left\{it\Gamma(1,n)\right\}\big]\cdot E\left[\exp\left\{it\sum_{l=2}^{p\enn}\Gamma(1,n)\right\}\right]\right|\\
    \leq&\varrho\p(X,s^\enn)\left\|\exp\left\{it\Gamma(1,n)\right\}\right\|_2
    \left\|\exp\left\{it\sum_{l=2}^{p\enn}\Gamma(1,n)\right\}\right\|_2\\
    \leq&\varrho\p(X,s^\enn).
\end{aligned}
\]
The last inequality is valid since the random variables are real (and so the exponential has modulus equal to 1).  This inequality suggests that if for each $c$, $1\leq c\leq p^\enn-1$, the terms
\[\prod_{l=1}^cE\left[\exp\left\{it\Gamma(l,n)\right\}\right]
E\left[\exp\left\{it\sum_{l=c+1}^{p\enn}\Gamma(l,n)\right\}\right]\]
are added and subtracted inside the absolute value brackets of \eqref{eq3.5.1}, then the triangle inequality implies that
\[
\left|\phi^\enn(t)-\psi^\enn(t)\right|\leq p^\enn \varrho\p(X,s^\enn)
\]
The definition of $p^\enn$ in \eqref{eq3.4.3}, however, implies that
\[p^\enn \varrho\p(X,s^\enn)\leq\sqrt{\varrho\p(X,s^\enn)}\to0\]
as $n\to\infty$.  This completes the proof of Lemma \ref{Lemma3.5}.\qed

\subsection{Application of Lyapounov's Condition to Reduction II}\label{C3S6}
It has just been shown by Lemma \ref{Lemma3.5} that it is sufficient to consider independent copies of the big blocks in order to analyze their end behavior.  In this section, Lyapounov's condition for the CLT will be applied to these independent blocks.\\

The Lyapounov condition for the $\DD$'s looks like the following:
\begin{lemma}\label{Lemma3.6}In the current context,
\begin{equation}\label{eq3.6.1}
    \lim_{n\to\infty}\frac{\sum_{l=1}^{p\enn}E|\DD|^4}{\left(\sum_{l=1}^{p\enn}E|\DD|^2\right)^2}=0.
\end{equation}
\end{lemma}

The proof of Lemma \ref{Lemma3.6} will commence in \S\ref{C3S6.5}:
\begin{claim}\label{Claim3.6}
With all the current notations, and with $\GG$ defined as in equation \eqref{Gamma},
\begin{equation}\label{eq3.6.2}
    \lim_{n\to\infty}\frac{\sum_{\smsp l}E\left|\Gamma(l,n)\right|^2}{\frac12f\lam\|\emph{\bv}\|^2}=1.
\end{equation}
\end{claim}

\noindent\textbf{Proof of Claim \ref{Claim3.6}:}  \\\\
Recall the definition of $\xi(l,n)$ in \eqref{xis}.  Also recall from \eqref{eq3.4.0b} that $p^\enn r^\enn/v_1^\enn\to1$.  It is therefore possible to interpret the sum
\[
\sum_{l=1}^{p\enn}E\Big(\Gamma(l,n)+\xi(l,n)\Big)^2
\]
as a type of Cesaro mean.  (To see this, you have to remember that implicit in the definition of \GG and $\xi(l,n)$ is a factor of $1/\sqrt{V^\enn}$, and then you can apply Miller's Lemma (Lemma \ref{LemmaM}) to the individual summands.)  Therefore,
\begin{equation}\label{eq3.6.3}
    \lim_{n\to\infty}\sum_{l=1}^{p\enn}E\Big(\Gamma(l,n)+\xi(l,n)\Big)^2=\frac12f\lam\|\bv\|^2.
\end{equation}

The trick now is to analyze the two quantities $\sum_lE\Big(\GG\Big)^2$ and $\sum_lE\Big(\xi(l,n)\Big)^2$.
Because of the Rosenthal Inequality (Theorem \ref{Brad2}),
\begin{equation}\label{eq3.6.4}
\begin{aligned}
E\Big(\GG\Big)^2\leq\frac{C}{V^\enn}\sum_{k\in\bB}
E\left[\left(\sum_{j=1}^ma_j\Re B_{k,q}^{\lambda\jenn}+b_j\Im B_{k,q}^{\lambda\jenn}\right)^2\right].
\end{aligned}
\end{equation}
Claim \ref{Claim2} now implies that
\begin{equation}\label{eq3.6.6}
\begin{aligned}
E\left[\left(\sum_{j=1}^ma_j\Re B_{k,q}^{\lambda\jenn}+b_j\Im B_{k,q}^{\lambda\jenn}\right)^2\right]\leq2m\|\bv\|^2\sum_{j=1}^m E \left|B_{k,q}^{\lambda\jenn}\right|^2
\end{aligned}
\end{equation}
By definition, $E \left|B_{k,q}^{\lambda\jenn}\right|^2<E|X_0|^2$. Putting this together with \eqref{eq3.6.4} and \eqref{eq3.6.6},
\begin{equation}\label{eq3.37}
\begin{aligned}
 \sum_{l=1}^{p\enn}E\Big(\GG\Big)^2\leq2Cm^2\|\bv\|^2\frac{p\enn\Card\big(\bB\big)}{V^\enn}E|X_0|^2.
\end{aligned}
\end{equation}
However, $\Card\big(\bB\big)=r^\enn\cdot v_2^\enn\cdots v_d^\enn$ and so since $r^\enn p^\enn/v_1^\enn\to1$
\[
\sup_n\sum_{l=1}^{p\enn}E\Big(\GG\Big)^2<\infty.
\]
What is more, \eqref{eq3.6.4} and \eqref{eq3.6.6} imply that 
\begin{equation}\label{eq3.41}
E\Big(\GG\Big)^2\leq 2Cm^2\|\bv\|^2 \frac{\Card\big(\bB\big)}{V^\enn}E|X_0|^2\sim \frac{C_1}{p\enn}
\end{equation}

With pretty much the same proof as \eqref{eq3.37}
\begin{equation}\label{eq3.38}
    \sum_{l=1}^{p\enn}E\Big(\xi(l,n)\Big)^2\leq\frac{2Cm\|\bv\|^2}{V^\enn}
    \sum_{l=1}^{p\enn}\sum_{\{k\in\bB\}}\sum_{j=1}^mE \left|T_{k,q}^{\lambda\jenn}\right|^2.
\end{equation}
Because of \eqref{eq3.2.1}, $E \left|T_{k,q}^{\lambda\jenn}\right|^2=E \left|T_{k,q}^{\lambda(1,n)}\right|^2$.  There is also a natural correspondence between the subscripts of $\bB$ and $\mathbb{B}(1,n)$; indices within each slice correspond if they are in the same position relative to their respective slices.  For example, each $k\in\bB$ corresponds to $k\p=k-((l-1)(s^\enn+r^\enn),0,\dots,0)\in\mathbb{B}(1,n)$.  Moreover, $\bk>\langle k\p\rangle$, (see \S3.2) hence
$E\left|X_{k\p}^{(\mu)}\mathbb{I}_{\{|X_{k\p}|>\langle k\p\rangle^q\}}\right|^2<E\left|X_k^{(\mu)}\mathbb{I}_{\{|X_k|>\bk^q\}}\right|^2$.
These facts, taken with the simple probabilistic inequality $Var(X)\leq E|X|^2$, imply that

\begin{equation}\label{eq3.39}
\begin{aligned}
    \sum_{l=1}^{p\enn}E\Big(\xi(l,n)\Big)^2&\leq\frac{2Cm^2\|\bv\|^2}{p^\enn r^\enn v_2^\enn\cdots v_d^\enn}
    \sum_{l=1}^{p\enn}\sum_{\{k\in\bB\}}E \left|T_{k,q}^{\lambda(1,n)}\right|^2\\
    &\leq \frac{2Cm^2\|\bv\|^2p^\enn}{p^\enn r^\enn v_2^\enn\cdots v_d^\enn}
  \sum_{\{k\in\mathbb{B}(1,n)\}}E\left|X_k^{(\mu)}\mathbb{I}_{\{|X_{k}|>\langle k\rangle^q\}}\right|^2.
\end{aligned}
\end{equation}
But, since $E|X_0|^2<\infty$ and $X$ is stationary, the summands in \eqref{eq3.39} satisfy the conditions of Lemma \ref{LemmaZ}.  Therefore,
\[
\sum_{l=1}^{p\enn}E\Big(\xi(l,n)\Big)^2\to0.
\]

Implicit in \eqref{eq3.39} is another uniform bound (that is, uniform with respect to the block number $l$) on the individual $\xi(l,n)$'s:
\begin{equation}\label{eq3.40}
   E\Big(\xi(l,n)\Big)^2\leq\frac{C_2}{r\enn p\enn v_2^\enn\cdots v_d^\enn}
  \sum_{\{k\in\mathbb{B}(1,n)\}}E\left|X_k^{(\mu)}\mathbb{I}_{\{|X_{k\p}|>\|k\p\|^q\}}\right|^2.
\end{equation} However, the summands satisfy the conditions of Lemma \ref{LemmaZ}, and so the right-hand side of \eqref{eq3.40} is $\sim \frac{C_2}{p\enn}\cdot o(1)$

The proof is nearly over; look again at \eqref{eq3.6.3}:
\[
\lim_{n\to\infty}\sum_{l=1}^{p\enn}E\Big(\Gamma(l,n)+\xi(l,n)\Big)^2=\frac12f\lam\|\bv\|^2.
\]
Obviously, $E\Big(\Gamma(l,n)+\xi(l,n)\Big)^2=E\Big(\Gamma(l,n)\Big)^2+E\Big(\xi(l,n)\Big)^2+2E\left[\GG\cdot\xi(l,n)\right]$.
The only things we need to worry about are the covariance terms. Referring to \eqref{eq3.40} and \eqref{eq3.41}, apply the Cauchy-Schwarz inequality:
\[
\begin{aligned}
\sum_l\Big|E\left[\GG\cdot\xi(l,n)\right]\Big|&\leq\sum_l\sqrt{E\left[\GG\right]^2}\cdot \sqrt{E\left[\xi(l,n)\right]^2}\\
&\leq\frac{1}{p^\enn}\sum_{l=1}^{p\enn} \sqrt{C_1}\sqrt{C_2\cdot o(1)}
\end{aligned}
\]
Claim \ref{Claim3.6} therefore holds, since this is essentially a Cesaro mean.  \qed
\subsection{Proof that Lyapounov's Condition Holds}\label{C3S6.5}

It is now possible to finish the proof of the lemma from the previous section.\\

\noindent\textbf{Proof of Lemma \ref{Lemma3.6}:}  A restatement of the desired result (Equation \eqref{eq3.6.1}) would be useful at this point:

\[
\lim_{n\to\infty}\frac{\sum_{l=1}^{p\enn}E|\DD|^4}{\left(\sum_{l=1}^{p\enn}E|\DD|^2\right)^2}=^?0.
\]
Claim \ref{Claim3.6} and the definition of $\DD$ imply that this is equivalent to proving
\[
\lim_{n\to\infty}\frac{\sum_{l=1}^{p\enn}E|\GG|^4}{\left(\frac12f\lam\|\bv\|^2\right)^2}=^?0.
\]

Notice first of all that

\begin{equation}\label{3.71}
\sum_{l=1}^{p\enn}E|\GG|^4\leq\frac{\displaystyle{\sum_{l=1}^{p\enn}}
E\left(\displaystyle{\sum_{k\in\bB}\sum_{j=1}^ma_j}\Re B_{k,q}^{\lambda\jenn}+b_j\Im B_{k,q}^{\lambda\jenn}\right)^4}
{\left(r\enn p\enn v_2^\enn\cdots v_d^\enn\right)^2}.
\end{equation}
This is a statement about moments, and in such situations, the Rosenthal inequality (Theorem 1.3.6)is very useful:
\begin{equation}\label{3.72}
\begin{aligned}
  &E\left({\sum_{k\in\bB}\sum_{j=1}^ma_j}\Re B_{k,q}^{\lambda\jenn}+b_j\Im B_{k,q}^{\lambda\jenn}\right)^4\\
  \leq& C\sum_{k\in\bB}E\left(\sum_{j=1}^ma_j\Re B_{k,q}^{\lambda\jenn}+b_j\Im B_{k,q}^{\lambda\jenn}\right)^4\\&\hspace{1in}+C\left(\sum_{k\in\bB}E\left(\sum_{j=1}^ma_j\Re B_{k,q}^{\lambda\jenn}+b_j\Im B_{k,q}^{\lambda\jenn}\right)^2\right)^2
\end{aligned}
\end{equation}
According to Claim \ref{Claim1},
\begin{equation}\label{3.73}
    \left(\sum_{j=1}^ma_j\Re B_{k,q}^{\lambda\jenn}+b_j\Im B_{k,q}^{\lambda\jenn}\right)^4
\leq m^2 \left(\sum_{j=1}^m(a_j^2+b_j^2)\left|B_{k,q}^{\lambda\jenn}\right|^2\right)^2.
\end{equation}
Recall that $\left|B_{k,q}^{\lambda\jenn}\right|\leq (k_1\cdot k_2\cdots k_d)^{q}$, and apply Claim \ref{Claim1} again:
\begin{equation}\label{3.74}\begin{aligned}
    m^2 \left(\sum_{j=1}^m(a_j^2+b_j^2)\left|B_{k,q}^{\lambda\jenn}\right|^2\right)^2
&\leq m^4\max_j\left\{a_j^4,b_j^4\right\}(k_1\cdot k_2\cdots k_d)^{4q}.\\
&\leq m^4\max_j\left\{a_j^4,b_j^4\right\}\left(V^\enn\right)^{4q}.
\end{aligned}\end{equation}
As shown before,
\begin{equation}\label{375}
\left(\sum_{j=1}^ma_j\Re B_{k,q}^{\lambda\jenn}+b_j\Im B_{k,q}^{\lambda\jenn}\right)^2
\leq m^2\max_j\left\{a_j^2,b_j^2\right\}E|X_0|^2.
\end{equation}
Combining \eqref{3.71}-\eqref{375}, it holds that
\[\label{376}
\begin{aligned}
&\sum_{l=1}^{p\enn}E|\GG|^4\\\leq&
\frac{\displaystyle{\sum_{l=1}^{p\enn}}
\left(\displaystyle{\sum_{k\in\bB}}C_1\left(V^\enn\right)^{4q}
+\left(\sum_{k\in\bB}C_2E|X_0|^2\right)^2\right)}
{\left(v_1^\enn v_2^\enn\cdots v_d^\enn\right)^2}\\
=&\frac{p\enn\left(C_1r\enn v_2^\enn\cdots v_d^\enn \cdot\left(V^\enn\right)^{4q}
+\left(C_2r\enn v_2^\enn\cdots v_d^\enn E|X_0|^2\right)^2\right)}{\left(v_1^\enn v_2^\enn\cdots v_d^\enn\right)^2}.
\end{aligned}
\]
If you now recall that $\lim_{n}(p\enn r\enn/v_1^\enn)=1$, and that $0<q<1/4$, it is not a problem to see that
\[
\lim_n\frac{p^\enn\left(C_1r\enn v_2^\enn\cdots v_d^\enn \cdot\left(V^\enn\right)^{4q}
\right)}{\left(p\enn r\enn v_2^\enn\cdots v_d^\enn\right)^2}=0.
\]
Moreover, since $p^\enn\to\infty$,
\[
\lim_n\frac{p\enn\left(C_2r\enn v_1^\enn\cdots v_d^\enn E|X_0|^2\right)^2}{\left(p\enn r\enn v_2^\enn\cdots v_d^\enn\right)^2}=0.
\]
This implies that $\lim_{n\to\infty}\frac{\sum_{l=1}^{p\enn}E|\GG|^4}{\left(\frac12f\lam\|\bv\|^2\right)^2}=0$, and hence that Lemma \ref{Lemma3.6} is true. \qed

\subsection{Proof of Reduction I and the Main Result}\label{C3S7}
In this section, the proofs of Theorem \ref{N3} and Lemma \ref{qLemma} will be completed.\\

\noindent{\textbf{Proof of Lemma \ref{qLemma}, Equation \eqref{Main3.2}:}} Lemma \ref{Lemma3.6}, and the Lyapounov CLT (see \cite{1}, Theorem 27.3) prove that
\[
\frac{\sum_{l=1}^{p\enn}\DD}{\sqrt{\frac12f\lam\|\bv\|^2}}\Rightarrow N(0,1).
\]
However, Lemma \ref{Lemma3.5} and \cite{1}, Theorem 26.3 imply that
\[
\frac{\sum_{l=1}^{p\enn}\GG}{\sqrt{\frac12f\lam\|\bv\|^2}}\Rightarrow N(0,1).
\]
Next, Lemma \ref{Lemma3.4} and Slutsky's Theorem prove that
\[
\frac{G\left(\bv,\bS_{n,q}^0\right)}
{\sqrt{\frac12f\lam\|\bv\|^2V^\enn}}\Rightarrow N(0,1).
\]
Since $\bv\in\R^k$ was arbitrary, the Cramer-Wold device implies that Lemma \ref{qLemma} holds.\qed\\

\noindent{\textbf{Conclusion of Proof of Theorem \ref{N3}:}  Since
\[G\left(\bv,\bS_n^0\right)
=G\left(\bv,\bS_{n,q}^0\right)+G\left(\bv,\left(Z_{n,q}^{\lambda(1,n)},\dots,Z_{n,q}^{\lambda(m,n)}\right)\right)\]
(it may be necessary to refer to the notations of \S\ref{C3S2}), Equation \eqref{Secondary3.2} of Lemma \ref{qLemma} and Slutsky's Theorem imply the main result. \qed
\newpage\smsp

\section{The Connection to Periodograms}\label{C4}

For the remainder of this paper, it is assumed implicitly that the sequence of vectors $v^\enn$ satisfies
\[\lim_{n\to\infty}\min\{v_1^\enn,v_2^\enn,\dots,v_d^\enn\}=\infty\]

The first two results are basically applications of Theorem \ref{N3} and the Mapping Theorem (Theorem \ref{mapp}).\\

The first result is trivial:

\begin{corollary}\label{CorA}
Let $d$ be a positive integer and suppose that $X$ is a $d$-dimensional CCSS random field such that $E|X_k|^2=\sigma<\infty$.  Let $f\lam:=f(e^{i\lambda})$ be the (continuous) spectral density of $X$. Let $\lambda\in\at$, and let $\{\lambda^\enn\}_{n=1}^\infty$ be a sequence which converges to $\lambda$.  Then
\[\frac1{\sqrt{V^\enn}}\left(\Re S_n^{\lambda\enn},\Im S_n^{\lambda\enn}\right)\Rightarrow N\left(0,\Upsilon_1^\lam\right).\]
\end{corollary}

Now consider the function $g:\R^{2m}\to\R^m$ defined by $g(x_1,x_2,\dots,x_{2m})=(x_1^2+x_2^2,\dots,x_{2m-1}^2+x_{2m}^2)$ it is continuous, and so its collection of discontinuities $D_g:=\{\textbf{x}\in\R^{2m}:g \text{ discontinuous at } \textbf{x}\}$ is empty.  Hence if $\nu$ (nu) is the measure defined by the distribution of \emph{any} multivariate normal, $\nu(D_g)=0$.  This suggests the following corollary:

\begin{corollary}\label{CorB}
Under the same hypotheses as Theorem \ref{N3},
\[
\left(I_n^{\lambda(1,n)},\dots,I_n^{\lambda(m,n)}\right)\Rightarrow \left(\chi_1^\lam,\chi_2^\lam,\dots,\chi_m^\lam\right)
\]
where the $\chi_j^\lam$ are independent exponential $(f\lam)$ random variables.
\end{corollary}

\pf A standard fact from probability theory is that if $N_1$ and $N_2$ are independent $N(0,\sigma)$ random variables, then $N_1^2+N_2^2$ has a scaled chi-squared distribution with two degrees of freedom, which is an exponential random variable with mean $2\sigma$.  Therefore, Theorem \ref{N3}, the Mapping Theorem, and the remarks immediately preceding the statement of Corollary \ref{CorB} prove the result.\qed\\

\rmk Notice the special case of Corollary \ref{CorB}: if $\{\lambda^\enn\}_{n=1}^\infty$ is a sequence which converges to $\lambda$, then $I_n^{\lambda\enn}\Rightarrow \chi^\lam$.\\

\end{document}